\documentclass{article}
\usepackage{latexsym}
\usepackage{epsfig}

\begin{document}

\title{\huge The Residual Intersection Formula of Type II Exceptional Curves}
\author{Ai-Ko Liu\footnote{Current Address: 
Mathematics Department of U.C. Berkeley}\footnote{
 HomePage:math.berkeley.edu/$\sim$akliu}}
\date{June, 16}
\maketitle
\newtheorem{main}{Main Theorem}
\newtheorem{theo}{Theorem}
\newtheorem{lemm}{Lemma}
\newtheorem{prop}{Proposition}
\newtheorem{rem}{Remark}
\newtheorem{cor}{Corollary}
\newtheorem{mem}{Examples}
\newtheorem{defin}{Definition}
\newtheorem{axiom}{Axiom}
\newtheorem{conj}{Conjecture}
\newtheorem{assum}{Assumption}
\medskip

\section{Preliminary}

 This paper is a part of the program [Liu1], [Liu3], [Liu4], [Liu5], [Liu6],
 [Liu7],
to understand the family Seiberg-Witten
theory and its relationship with the enumeration of nodal or singular curves
 in linear systems of algebraic surfaces. In [Liu1] a symplectic approach
 to the universality theorem is given. In [Liu6] the algebraic geometric
approach is given. In [Liu7] this result has been interpreted as an enumerative
 Riemann-Roch formula probing the non-linear information of the linear systems.

The universality theorem implies that for $5n-1$ very ample line bundle
$L\mapsto M$, the ``number of $n$-node nodal curves'' in a generic
 $n$ dimensional linear sub-system of $|L|$ can be expressed as a
universal polynomial of the characteristic numbers $c_1^2(L), c_1({\bf K}_M)
\cdot c_1(L)$, $c_1^2({\bf K}_M)$ and $c_2(M)$, in the spirit of the
 surface Riemann-Roch formula. On the other hand, for $L$ not sufficiently
 very ample, the actual virtual number of nodal curves differs from the
 universal formula predicted by G$\ddot{o}$ttsche [Got]. In [Liu2] the corrections
 from the type $II$ exceptional classes have been interpreted as a non-linear
 analogue of second sheaf cohomology.

 In this paper, we build up the theory of type $II$ exceptional classes,
 parallel to the type $I$ theory built up in [Liu1], [Liu5] and [Liu6].

 One major application of the type $II$ theory is to define the
 ``virtual number of nodal curves'' in $|L|$ on algebraic surfaces 
without any condition on $L$. 
 
A direct application of our theory is to argue the vanishing result of 
 type $II$ contributions on universal families of
 $K3$ surfaces. Once this is achieved, the ``virtual numbers of nodal curves''
 on $K3$ are equal to the polynomials constructed from
 the universality theorem [Liu6].

Another interesting application of the theory of type $II$ exceptional 
classes to enumerative problem is the solution of Harvey-Moore conjecture [Liu2]
 on the enumeration of nodal curves on Calabi-Yau K3 fibrations.

  The layout of the paper is as the following. 
 In section \ref{section; type2}, we review the algebraic family Kuranishi 
models of type $II$ exceptional classes.

 Then in section \ref{section; kura}, we construct the Kuranishi models explicitly.
 In section \ref{section; blowup}, we consider the blowup construction of
the algebraic family Seiberg-Witten invariants and prove the main theorem 
of the paper on the mixed algebraic family Seiberg-Witten invariants 
attached to a finite collection of type $II$ exceptional classes. 

 The following is an abbreviated form of our
 main theorem of the paper, stated in a less technical term.
Please refer to theorem \ref{theo; degenerate} on 
page \pageref{theo; degenerate} for the more 
complete statement.

\begin{main}\label{main; 1}
Given an algebraic family of algebraic surfaces $\pi:{\cal X}\mapsto B$ and a
 finite collection of type $II$ exceptional 
\footnote{For the definition of exceptional classes, please
 consult definition \ref{defin; exception} on page 
\pageref{defin; exception}.} classes, $e_{II; 1}$, 
$e_{II; 2}$, $e_{II; 3}$, $\cdots$, $e_{II; p}$ along ${\cal X}\mapsto B$ 
satisfying $e_{II; i}\cdot
 e_{II; j}\geq 0$ for $i\not=j$, then the localized (excess) contribution of
 the algebraic family Seiberg-Witten invariant  
${\cal AFSW}_{{\cal X}\mapsto B}(1, C)$ along the locus of co-existence 
 $\times_B^{1\leq i\leq p}\pi_i(\times_B^{1\leq i\leq p}{\cal M}_{e_{II; i}})$ 
  of type $II$ exceptional classes
is well defined as a mixed algebraic family Seiberg-Witten invariant 
 of $C-\sum_{1\leq i\leq p}e_{II; i}$, manifestly a topological invariant
 independent of the
 choices of the family Kuranishi models and the possible deformations
 of the family $\pi:{\cal X}\mapsto B$. As a direct consequence,  
the residual contribution of
 the family invariant ${\cal AFSW}_{{\cal X}\mapsto B}(1, C)$, which is
 ${\cal AFSW}_{{\cal X}\mapsto B}(1, C)$ subtracted by the localized excess
contribution, is well defined. 
\end{main}

The above theorem also works for the mixed
 invariants ${\cal AFSW}_{{\cal X}\mapsto 
B}(\eta, C)$, $\eta\in {\cal A}_{\cdot}(B)$.

The above theorem generalizes the theory of type $I$ exceptional curves developed 
in [Liu5] and [Liu6] to the cases when the moduli spaces of exceptional curves
 are not regular.

\medskip

 At the end of the paper, we outline the procedure to extend the scheme
 to an inductive scheme on hierarchies of finite collections of 
type $II$ classes. Then we apply the inductive scheme
 to the universal families of $K3$ surfaces and argue the vanishing results
 on $K3$ universal families.

\medskip

\section{The Review of Algebraic Family Kuranishi Models of Type II Exceptional
 Classes}\label{section; type2}

\bigskip

  Recall that a type $I$ exceptional class $e_i=E_i-\sum_{j_i}E_{j_i}$ 
of the fiber bundle of universal spaces
 $M_{n+1}\mapsto M_n$ has the following two crucial properties:

\bigskip

\noindent 1. The family moduli space of $e_i$ is smooth of codimension
 $-{e_i^2-c_1({\bf K}_{M_{n+1}/M_n})\cdot e_i\over 2}$ $=-e_i^2-1$
 in $M_{n}$ and can be
 identified with the closure of an admissible stratum, 
$Y(\Gamma_{e_i})$, for a fan-like
 admissible graph $\Gamma_{e_i}$.

\begin{figure}
\centerline{\epsfig{file=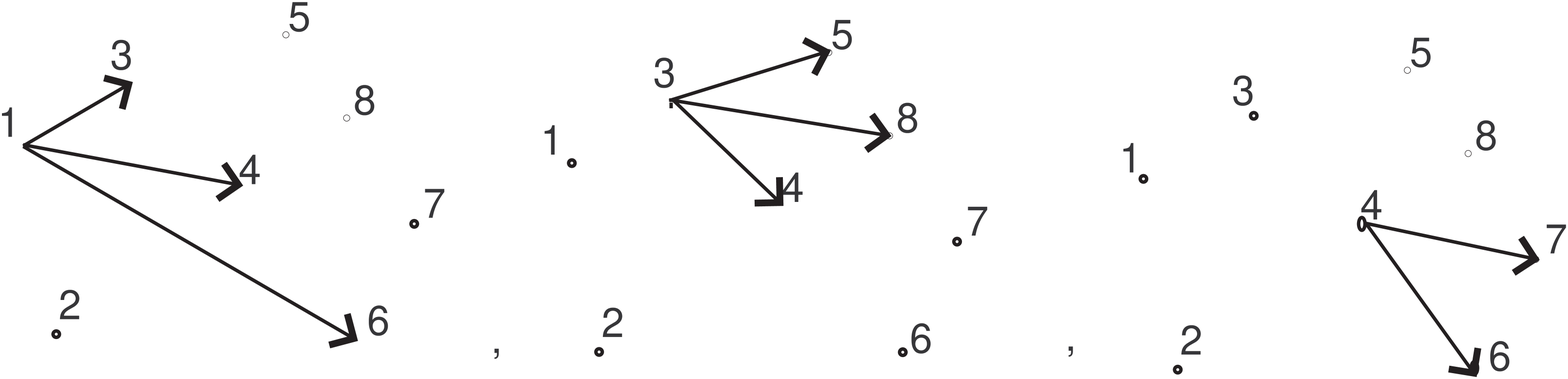,height=4cm}}
\centerline{fig.1}
\centerline{A few fan-like admissible graphs with eight vertexes}
\end{figure}
\label{fanlike}

\medskip

\noindent 2. Over each $b\in Y(\Gamma_{e_i})$, the class $e_i$ 
is represented by a
 unique holomorphic curve, called type $I$ exceptional curve,
 in the fiber $M_{n+1}|_b=M_{n+1}\times_{M_n}\{b\}$. Over a 
Zariski open subset
 $Y_{\Gamma_{e_i}}$ of 
$Y(\Gamma_{e_i})$, the curves representing $e_i$ are irreducible.
 Over a finite union of codimension one loci in 
$Y(\Gamma_{e_i})-Y_{\Gamma_{e_i}}$, the
 curves are disjoint unions of irreducible components, each of them 
 being an irreducible type $I$ exceptional curve.

\medskip 

 The fact that the fibrations of type $I$ curves are smooth and ``universal''
 has played a crucial role in understanding the universal nature of the
 ``universality'' theorem [Liu1], [Liu6], as a natural extension of surface
 Riemann-Roch theorem in enumerative geometry.

 The above properties about type $I$ exceptional curves
 are the consequence of the existence of a ``canonical'' algebraic
 family Kuranishi model of $e_i$.
 These two properties have been used extensively in the 
proofs of the universality theorem [Liu1], [Liu6]. In this paper, we
 develop the necessary
 algebraic technique to deal with type $II$ exceptional curves.

\bigskip

 Recall the following definition of exceptional classes which had 
already appeared in
 [Liu4], [Liu7],

\begin{defin}\label{defin; exception}
 A class $e$  
is said to be exceptional over $\pi:{\cal X}\mapsto B$ 
if it satisfies the following
 conditions:

\noindent (i). The fiberwise self-intersection number $\int_B\pi_{\ast}(e^2)<0$.

\medskip

\noindent (ii). The relative degree $deg_{{\cal X}/B}e>0$.

\end{defin}

\medskip

\begin{defin}\label{defin; type2}
 Consider the universal family ${\cal X}=M_{n+1}\mapsto B=M_n$. An exceptional
 class $e$ is said to be a type $II$ exceptional class if 
 $e$ does not lie in the kernel of ${\cal A}_{\cdot}(M_{n+1})\mapsto
 {\cal A}_{\cdot}(M\times M_n)$.
\end{defin}

 For a general (non-universal) fiber bundle ${\cal X}\mapsto B$,
 we also use the term ``type $II$ exceptional class'' in referring exceptional
 classes of the fibration ${\cal X}\mapsto B$.
 In this paper, we often use $e_{II}$ with a subscript $II$
 to denote a type $II$ exceptional class.

To simplify our discussion and get to the key point, we impose an 
additional condition on $e_{II}$. 

\begin{assum}( Assumption on $e_{II}$)\label{assum; 2}
 $deg_{{\cal X}/B}(c_1({\bf K}_{{\cal X}/B})-e_{II})<0$.
\end{assum}

 This implies that ${\cal R}^2\pi_{\ast}\bigl({\cal E}_{e_{II}}\bigr)=0$ by
 relative Serre duality.

\medskip

For a type $II$ exceptional class $e_{II}$ of the fiber bundle
 ${\cal X}\mapsto B$, there is usually no canonical choices of the
 algebraic family Kuranishi models. We know that when
 a curve representing $e_{II}$ is irreducible, it must be unique in the same
 fiber. However the reducible representatives may even contain some irreducible
 component with a non-negative self-intersection. Thus,
 the following general symptoms have to be kept in mind,

\bigskip \label{bizzare}

\noindent 1'. The family moduli space of $e_{II}$, 
${\cal M}_{e_{II}}\mapsto B$ (discussed
 in more detail below), may be not of the expected algebraic family
 dimension $dim_{\bf C}B+p_g+{e_{II}^2-c_1({\bf K})\cdot e_{II}\over 2}$. 
 The sub-locus $\subset 
{\cal M}_{e_{II}}$ over which the universal curve representing 
 $e_{II}$ remains irreducible may not be open and dense in ${\cal M}_{e_{II}}$
 and can be even empty sometimes.

\medskip

\noindent 2'. The projection map ${\cal M}_{e_{II}}\mapsto B$
 is usually not a closed 
immersion. And for each $b\in B$, the fiber of ${\cal M}_{e_{II}}\mapsto B$
 above $b$, ${\cal M}_{e_{II}}|_b$, parametrizing all 
the curves dual to $e_{II}$ in the
 fiber ${\cal X}_b$, can be of positive dimension. This means that
 there can be (uncountably) many representatives of $e_{II}$
 above this given $b\in B$. By the exceptionality condition on $e_{II}$, 
this can occur only when the representative contains more than
 one irreducible component.

\medskip

  Our task is to develop a version of residual intersection formula of
 the algebraic family Seiberg-Witten invariant based on the general symptoms 
1' and 2'. We will demonstrate that there is a well-defined theory of 
 type $II$ exceptional classes parallel
 to the theory of type $I$ exceptional classes.  The basic tool we will use
 is the construction of the algebraic family Kuranishi models of $e_{II}$.

\subsection{\bf A Short Review About the Kuranishi-Model of $e_{II}$}
\label{subsection; review}

\bigskip

 In the following, we give a short review about the construction of the
algebraic family Kuranishi models of $e_{II}$ and discuss their basic 
properties. Because our main application is about the universal families,
 we assume that ${\cal R}^2\pi_{\ast}{\cal O}_{\cal X}$ is isomorphic to
 ${\cal O}_B^{p_g}$. In this case\footnote{For the
definition of the formal excess base dimension $febd$, please
 consult definition 4.5. of [Liu3] for more details.}, 
$febd(e_{II}, {\cal X}/B)=p_g$. We will 
assume implicitly
 in most of the current paper that $febd(e_{II}, {\cal X}/B)=p_g$ 
to simplify our discussion.

 Because in each fiber algebraic surface 
of the fiber bundle ${\cal X}\mapsto B$, 
the curve (divisor) representing $e_{II}$ may not be unique, 
 we consider the base 
space, ${\cal T}_B({\cal X})$, of relative $Pic^0$ tori parametrizing all the
 holomorphic structures of ${\cal E}_{e_{II}}$ over $B$. 
 Let $({\bf V}_{II}, {\bf W}_{II}, \Phi_{{\bf V}_{II}{\bf W}_{II}})$
 be an algebraic family Kuranishi model of $e_{II}$, defined over
 ${\cal T}_B({\cal X})$ and let $\Phi_{{\cal V}_{II}{\cal W}_{II}}:{\cal V}_{II}
\mapsto {\cal W}_{II}$ be the corresponding morphism of locally 
free sheaves. Then by its defining property and the simplifying 
assumption \ref{assum; 2} on $e_{II}$ (on page \pageref{assum; 2}), we have 
$Ker(\Phi_{{\cal V}_{II}{\cal W}_{II}})\cong {\cal R}^0\pi_{\ast}{\cal E}_{e_{II}}$
 and $Coker(\Phi_{{\cal V}_{II}{\cal W}_{II}})\cong 
{\cal R}^1\pi_{\ast}{\cal E}_{e_{II}}$, where
 ${\cal R}^i\pi_{\ast}\bigl({\cal E}_{e_{II}}\bigr)$ over 
${\cal T}_B({\cal X})$ is the $i$-th derived
 image sheaf of ${\cal E}_{e_{II}}$ over ${\cal X}\times_B {\cal T}_B({\cal X})$
 along $\pi:{\cal X}\times_B {\cal T}_B({\cal X})\mapsto {\cal T}_B({\cal X})$.
 
   The kernel $Ker(\Phi_{{\bf V}_{II}{\bf W}_{II}})$ defines an algebraic sub-cone 
${\bf C}_{e_{II}}$ over
 ${\cal T}_B({\cal X})$, and its projectification ${\bf P}({\bf C}_{e_{II}})$
 is nothing but the algebraic family moduli space ${\cal M}_{e_{II}}$ 
of $e_{II}$ over $B$.

 The bundle map ${\bf V}_{II}\mapsto {\bf W}_{II}$ induces a global
 section $s_{II}$ 
of ${\bf H}\otimes \pi_{{\bf P}({\bf V}_{II})}^{\ast}{\bf W}_{II}$ over
 ${\bf P}({\bf V}_{II})$
 such that ${\cal M}_{e_{II}}$ can be identified with the zero locus
 $Z(s_{II})$ of the section $s_{II}\in \Gamma({\bf P}({\bf V}_{II}), 
 {\bf H}\otimes \pi_{{\bf P}({\bf V}_{II})}^{\ast}{\bf W}_{II})$. 

  As $Z(s_{II})$ may not be regular or be of the expected dimension, we
 have to rely on intersection theory [F] to construct the virtual fundamental
class of ${\cal M}_{e_{II}}$.

By applying the concept of localized top Chern class
 on page 244, section 14.1 of [F], $$[{\cal M}_{e_{II}}]_{vir}\doteq 
{\bf Z}(s_{II})\in 
{\cal A}_{dim_{\bf C}{\cal T}_B({\cal X})+rank_{\bf C}{\bf V}_{II}-1-rank_{\bf C}
{\bf W}_{II}}({\cal M}_{e_{II}})$$ defines a unique cycle class representing
 the ``virtual fundamental class'' of the family moduli space ${\cal M}_{e_{II}}$
 graded by the expected algebraic family Seiberg-Witten dimension
\footnote{This works for $e_{II}$ with $febd(e_{II}, {\cal X}/B)=p_g$.},

$$\hskip -.2in
ed=dim_{\bf C}{\cal T}_B({\cal X})+rank_{\bf C}{\bf V}_{II}-1-rank_{\bf C}
{\bf W}_{II}=dim_{\bf C}B+q+(p_g-q+{e_{II}^2-e_{II}\cdot c_1({\bf K}_{{\cal X}/B})
\over 2})$$
$$=dim_{\bf C}B+p_g+{e_{II}^2-e_{II}\cdot c_1({\bf K}_{{\cal X}/B})
\over 2}.$$

 This virtual fundamental class, i.e. the localized top Chern class
 ${\bf Z}(s_{II})$, can be pushed-forward and mapped into
a global cycle class in 

\noindent 
${\cal A}_{dim_{\bf C}{\cal T}_B({\cal X})+rank_{\bf C}{\bf V}_{II}-1-rank_{\bf C}
{\bf W}_{II}}({\bf P}({\bf V}_{II}))$ induced 
by the inclusion ${\cal M}_{e_{II}}\mapsto
 {\bf P}({\bf V}_{II})$.
 
 As we expect, the virtual fundamental class ${\bf Z}(s_{II})$ localized in
 ${\cal A}_{\cdot}({\cal M}_{e_{II}})$ is
 independent to the choices of the algebraic family Kuranishi models of $e_{II}$.

\begin{prop}\label{prop; independent}
The cycle class ${\bf Z}(s_{II})\in {\cal A}_{\cdot}({\cal M}_{e_{II}})$ 
is independent to the
 choices of the algebraic family Kuranishi model 
$({\bf V}_{II}, {\bf W}_{II}, \Phi_{{\bf V}_{II}{\bf W}_{II}})$ of $e_{II}$.
\end{prop}

\noindent Proof: Consider the fiber square 

\[
\begin{array}{ccc}
 Z(s_{II}) & \longrightarrow &  {\bf P}({\bf V}_{II}) \\
 \Big\downarrow &  & \Big\downarrow \vcenter{%
\rlap{$\scriptstyle{\mathrm{s_{II}}}\,$}}\\
 {\bf P}({\bf V}_{II}) &  
\stackrel{s_{\pi_{{\bf P}({\bf V}_{II})}^{\ast}{\bf W}_{II}\otimes
 {\bf H}}}{\longrightarrow} & \pi_{{\bf P}({\bf V}_{II})}^{\ast}{\bf W}_{II}\otimes
 {\bf H}
\end{array}
\]

, where $s_{\pi_{{\bf P}({\bf V}_{II})}^{\ast}{\bf W}_{II}\otimes
 {\bf H}}$ is the zero cross section of 
$\pi_{{\bf P}({\bf V}_{II})}^{\ast}{\bf W}_{II}\otimes
 {\bf H}$. By the fact that ${\bf Z}(s_{II})=
s_{\pi_{{\bf P}({\bf V}_{II})}^{\ast}{\bf W}_{II}\otimes {\bf H}}^{!}
[{\bf P}({\bf V}_{II})]$, ${\bf Z}(s_{II})$ is nothing but the following 
localized contribution of
 top Chern class along $Z(s_{II})$, 

 $$\hskip -.2in
\{c_{total}(\pi_{{\bf P}({\bf V}_{II})}^{\ast}{\bf W}_{II}\otimes 
{\bf H}|_{Z(s_{II})})\cap 
 s_{total}(Z(s_{II}), {\bf P}({\bf V}_{II}))
\}_{dim_{\bf C}{\cal T}_B({\cal X})+rank_{\bf C}{\bf V}_{II}-1-rank_{\bf C}
{\bf W}_{II}}.$$

 It suffices to show that the above localized contribution of top Chern class
 has been independent of the choices of ${\bf V}_{II}$ and  ${\bf W}_{II}$.

\begin{lemm}\label{lemm; stable}
the localized contribution of top Chern class of 
$s:\pi_{{\bf P}({\bf E})}^{\ast}{\bf F}\otimes {\bf H}$ induced by
 $\sigma:{\bf E}\mapsto {\bf F}$ is invariant under the stabilization
 $\sigma\sim
\sigma\oplus id_{\bf G}:{\bf E}\oplus {\bf G}\mapsto {\bf F}\oplus {\bf G}$.
\end{lemm}

 The lemma is similar to lemma 5.3. in [Liu3].

\noindent Proof of the lemma: Under the smooth 
embedding ${\bf P}({\bf E})\hookrightarrow
 {\bf P}({\bf E}\oplus {\bf G})$, the normal bundle of ${\bf P}({\bf E})$ in
 ${\bf P}({\bf E}\oplus {\bf G})$ is isomorphic to the bundle 
$\pi_{{\bf P}({\bf E})}^{\ast}{\bf G}\otimes {\bf H}$, as 
 ${\bf P}({\bf E})$ can be viewed as the zero locus of a regular section
 of $\pi_{{\bf P}({\bf E}\oplus {\bf G})}^{\ast}{\bf G}\otimes {\bf H}$ 
induced by the bundle projection ${\bf E}\oplus {\bf G}\mapsto {\bf G}$. So
 the total Segre class 

$$s_{total}({\bf P}({\bf E}), {\bf P}({\bf E}\oplus {\bf G}))=
s_{total}(\pi_{{\bf P}({\bf E})}^{\ast}{\bf G}\otimes {\bf H}),$$

and

$$s_{total}(Z(s_{II}), {\bf P}({\bf E}\oplus {\bf G}))
=s_{total}({\bf P}({\bf E}),
 {\bf P}({\bf E}\oplus {\bf G}))\cap 
s_{total}(Z(s_{II}), {\bf P}({\bf E}))$$

$$=c_{total}(-\pi_{{\bf P}({\bf E})}^{\ast}{\bf G}\otimes {\bf H})\cap
 s_{total}(Z(s_{II}), {\bf P}({\bf E})).$$

 Thus 

$$
c_{total}(\pi_{{\bf P}({\bf E}\oplus {\bf G})}^{\ast}({\bf E}\oplus 
{\bf G})\otimes {\bf H}|_{Z(s_{II})})
\cap s_{total}(Z(s_{II}), {\bf P}({\bf E}\oplus {\bf G}))$$
$$=c_{total}(\bigl(\pi_{{\bf P}({\bf E})}^{\ast}({\bf E}\oplus 
{\bf G})\otimes {\bf H}-
\pi_{{\bf P}({\bf E})}^{\ast}{\bf G}\otimes {\bf H}\bigr)|_{Z(s_{II})})\cap 
 s_{total}(Z(s_{II}), {\bf P}({\bf E}))$$

$$=c_{total}(\pi_{{\bf P}({\bf E})}^{\ast}({\bf E}\otimes 
{\bf H}|_{Z(s_{II})})\cap 
 s_{total}(Z(s_{II}), {\bf P}({\bf E})).$$

 So the localized contribution of top Chern class is invariant under the
 stabilization. $\Box$

 Once the lemma is proved, we may show that the localized top Chern
 classes defined by any
 two algebraic family Kuranishi models $(\Phi_{{\bf V}_{II}{\bf W}_{II}}, 
{\bf V}_{II}, {\bf W}_{II})$ and $(\Phi_{{\bf V}_{II}'{\bf W}_{II}'}, 
{\bf V}_{II}', {\bf W}_{II}')$ are equal.

 In fact one may stabilize $(\Phi_{{\bf V}_{II}{\bf W}_{II}}, 
{\bf V}_{II}, {\bf W}_{II})$ into $(\Phi_{{\bf V}_{II}{\bf W}_{II}}\oplus
 id_{{\bf V}_{II}'}, 
{\bf V}_{II}\oplus {\bf V}_{II}', {\bf W}_{II}\oplus {\bf V}_{II}')$
 and $(\Phi_{{\bf V}_{II}'{\bf W}_{II}'}, 
{\bf V}_{II}', {\bf W}_{II}')$ into 
$(id_{{\bf V}_{II}}\oplus\Phi_{{\bf V}_{II}'{\bf W}_{II}'}, 
{\bf V}_{II}\oplus {\bf V}_{II}', {\bf V}_{II}\oplus {\bf W}_{II}')$, respectively,
  by applying lemma \ref{lemm; stable}. We find that the
localized top Chern classes are stabilized into 

$$\{c_{total}(\pi_{{\bf P}({\bf V}_{II}\oplus {\bf V}_{II}')}^{\ast}
({\bf W}_{II}\oplus {\bf V}_{II}')\otimes {\bf H})\cap
 s_{total}(i_1({\cal M}_{e_{II}}, {\bf P}({\bf V}_{II}\oplus {\bf V}_{II}'))\}_{ed}$$

 and 

$$\{c_{total}(\pi_{{\bf P}({\bf V}_{II}\oplus {\bf V}_{II}')}^{\ast}
({\bf W}_{II}'\oplus {\bf V}_{II})\otimes {\bf H})\cap
 s_{total}(i_2({\cal M}_{e_{II}}, {\bf P}({\bf V}_{II}\oplus {\bf V}_{II}'))\}_{ed},$$

respectively. Over here $i_1, i_2:{\cal M}_{e_{II}}\hookrightarrow 
{\bf P}({\bf V}_{II}\oplus {\bf V}_{II}')$ denote two different imbeddings
${\cal M}_{e_{II}}\subset {\bf P}({\bf V}_{II})\subset 
 {\bf P}({\bf V}_{II}\oplus {\bf V}_{II}')$ and 
${\cal M}_{e_{II}}\subset {\bf P}({\bf V}_{II}')\subset 
 {\bf P}({\bf V}_{II}\oplus {\bf V}_{II}')$, respectively.

Firstly, because both ${\cal V}_{II}-{\cal W}_{II}$ and 
${\cal V}_{II}'-{\cal W}_{II}'$ are equal to 
 ${\cal R}^0\pi_{\ast}\bigl({\cal O}_{\cal X}(e_{II})\bigr)-
{\cal R}^1\pi_{\ast}\bigl({\cal O}_{\cal X}(e_{II})\bigr)$ in 
$K_0({\cal T}_B({\cal X}))$, we have ${\bf W}_{II}\oplus {\bf V}_{II}'\equiv 
 {\bf W}_{II}'\oplus {\bf V}_{II}$ and the corresponding total Chern classes
 are equal.

 Secondly, to show that 
$\rho_1=s_{total}(i_1({\cal M}_{e_{II}}), {\bf P}({\bf V}_{II}\oplus {\bf V}_{II}'))$
 and $\rho_2=s_{total}(i_2({\cal M}_{e_{II}}), {\bf P}({\bf V}_{II}\oplus {\bf V}_{II}'))$ 
 are equal, we notice that $i_1, i_2$ are within a ${\bf P}^1$ pencil
 of imbeddings of ${\cal M}_{e_{II}}$ induced by 
 $j_{a, b}:{\bf C}_{e_{II}}\mapsto {\bf V}_{II}\oplus {\bf V}_{II}'$;
 $j_{a, b}(v)=aj_1(v)\oplus bj_2(v)$ for $v\in {\bf C}_{e_{II}}$, $(a, b)\in
 {\bf C}^2-(0, 0)$. 
Here $j_1:{\bf C}_{e_{II}}\mapsto
 {\bf V}_{II}$ and $j_2:{\bf C}_{e_{II}}\mapsto
 {\bf V}_{II}'$ are the imbeddings of abelian cones projectified into
 ${\cal M}_{e_{II}}\subset {\bf P}({\bf V}_{II}), {\bf P}({\bf V}_{II}')$, 
respectively.

 So we may consider the embedding ${\bf P}^1\times {\cal M}_{e_{II}}\mapsto 
 {\bf P}^1\times {\bf P}({\bf V}_{II}\oplus {\bf V}_{II}')$ and the total 
Segre class
 $s_{total}({\bf P}^1\times {\cal M}_{e_{II}}, 
 {\bf P}^1\times {\bf P}({\bf V}_{II}\oplus {\bf V}_{II}'))$ of its normal
 cone.

 It is clear that if we restrict the total 
Segre class to different $\{t\}\times {\cal M}_{e_{II}}$,
 $t\in {\bf P}^1$, the resulting class $\in {\cal A}_{\cdot}({\cal M}_{e_{II}})$
 is independent to $t$. When $t=0$ and $t=\infty$, we recover
$\rho_1$ and $\rho_2$, respectively. Thus $\rho_1=\rho_2$.

 Because the Chern classes and Segre classes are identified, so are the
 corresponding localized top Chern classes of ${\cal M}_{e_{II}}$.
 $\Box$ 

 \bigskip

Because the localized top Chern class is canonically defined, we will denote
 them by $[{\cal M}_{e_{II}}]_{vir}$.

 If we push-forward the cycle class ${\bf Z}(s_{II})=[{\cal M}_{e_{II}}]_{vir}$
 into
 ${\bf P}({\bf V}_{II})$, then by proposition 14.1(a) on page 244 of [F],
 the image cycle class 
is equal to the global $c_{top}(\pi_{{\bf P}({\bf V}_{II})}^{\ast}
{\bf W}_{II}\otimes {\bf H})\cap [{\bf P}({\bf V}_{II})]$.

\bigskip

\section{The Construction of Algebraic Family 
Kuranishi Models}\label{section; kura}

\bigskip

 In the previous section, we have discussed the independence of 
$[{\cal M}_{e_{II}}]_{vir}$ to the choices of the algebraic family
 Kuranishi models of $e_{II}$. In the following, we will first review
 the construction of the family Kuranishi models. When we perform the
 blowup/residual intersection theory
 construction of algebraic family Seiberg-Witten invariants in
 subsection \ref{subsection; stable}, these
explicitly constructed algebraic family Kuranishi models will play a crucial
role.

 As was mentioned earlier, we focus mostly on the $febd(e_{II}, {\cal X}/B)=p_g$
case in the following discussion.

\bigskip
\subsection{The Construction of family Kuranishi Model of the Type $II$
 Exceptional Curves}\label{subsection; type2K}

 In this subsection, we review the explicit construction of the algebraic 
family Kuranishi model
 of a type $II$ exceptional class $e_{II}$.

 Let $\pi:{\cal X}\mapsto B$ be a fiber bundle of algebraic surfaces and let
 ${\cal T}_B({\cal X})$ be the fiber bundle of the relative $Pic^0$ tori. 
As usual, let
 ${\cal E}_{e_{II}}\mapsto {\cal T}_B({\cal X})$ 
be the invertible sheaf corresponding to $e_{II}$. 

\medskip

 Let $D\subset
 {\cal X}$ be an ample effective divisor on 
${\cal X}$ and let $n$ be a sufficiently
 large integer.

\begin{lemm}\label{lemm; ample}
 Suppose that $|D|$ is chosen to be sufficiently very ample, 
 then the divisor $D$ in $|D|$ can be chosen such that
 the composition map $D\subset {\cal X}\mapsto B$ is of relative dimension one.
\end{lemm}

\noindent Proof of lemma \ref{lemm; ample}:
 For all the closed points $b\in B$, consider the ${\cal O}_{{\cal X}}(D)$-twisted 
short exact sequence,
 $$0\mapsto {\cal I}_{{\cal X}_b}(D)\mapsto {\cal O}_{\cal X}(D)\mapsto 
{\cal O}_{{\cal X}_b}(D)\mapsto 0.$$

 By theorem 1.5. of [Ko], we may replace $D$ by a suitably large 
multiple such that ${\cal R}^i\pi_{\ast}{\cal I}_{{\cal X}_b}(D)=0$,
 ${\cal R}^i\pi_{\ast}{\cal O}_{{\cal X}_b}(D)=0$, for 
 $i>0$. 

 So the derived exact sequence from the above short exact sequence 
generates
 a short exact sequence \footnote{Here $k(b)$ is the residue field
 of $b$.} for each $b\in B$,

$$0\mapsto H^0({\cal X}, {\cal I}_{{\cal X}_b}(D))\otimes k(b)
\mapsto H^0({\cal X}, {\cal O}_{\cal X}(D))\mapsto 
H^0({\cal X}_b, {\cal O}_{{\cal X}_b}(D))\mapsto 0.$$

  The space $H^0({\cal X}, {\cal I}_{{\cal X}_b}(D))\otimes k(b)$ 
 is the subspace of the global sections 
 $H^0({\cal X}, {\cal O}_{\cal X}(D))$ which restricts to the trivial
section to ${\cal X}_b$. When $b$ moves these vector spaces form
 a vector bundle, denoted by ${\bf U}$.
 Its rank can be calculated to be
 $$\chi({\cal X}, {\cal O}_{\cal X}(D))-\chi({\cal X}_b, {\cal O}_{{\cal X}_b}(D))
\ll dim_{\bf C}|D|-dim_{\bf C}B,$$

 if $D$ is chosen such that $\chi({\cal X}_b, {\cal O}_{{\cal X}_b}(D))\gg
 dim_{\bf C}B$.  If such an inequality has been achieved, then $dim_{\bf C}|D|$
 is much larger than the dimension of the projective space bundle 
${\bf P}({\bf U})$ over $B$.
 By choosing an element 
of $|D|-Im({\bf P}({\bf U}))$, it gives rise to a
 cross section which restricts to non-trivial
 sections to each fiber ${\cal X}_b$. So after replacing the original 
$D$ by the defining divisor $D$ of the
 chosen section,  the newly chosen $D$ intersects
 each ${\cal X}_b$ properly and cuts it down into a curve.  
The lemma is proved. $\Box$

\medskip
 Once we have such a carefully chosen $D$, we are ready to construct
 Kuranishi models of $e_{II}$.

 The following short exact sequence

$$\hskip -.3in
0\mapsto {\cal O}_{\cal X}\otimes {\cal E}_{e_{II}}\mapsto
 {\cal O}_{\cal X}(nD)\otimes {\cal E}_{e_{II}}\mapsto
{\cal O}_{nD}(nD)\otimes {\cal E}_{e_{II}}\mapsto 0$$

 implies

$$
0\mapsto {\cal R}^0\pi_{\ast}\bigl({\cal O}_{\cal X}\otimes {\cal E}_{e_{II}}\bigr)
\mapsto {\cal R}^0\pi_{\ast}\bigl({\cal O}_{\cal X}(nD)\otimes {\cal E}_{e_{II}}
\bigr)\mapsto 
{\cal R}^0\pi_{\ast}\bigl({\cal O}_{nD}(nD)\otimes {\cal E}_{e_{II}}\bigr)$$

$$\hskip -.4in \mapsto
{\cal R}^1\pi_{\ast}\bigl({\cal O}_{\cal X}\otimes {\cal E}_{e_{II}}
\bigr)\mapsto 0,$$

 and ${\cal R}^1\pi_{\ast}\bigl({\cal O}_{nD}(nD)\otimes {\cal E}_{e_{II}}\bigr)
\cong {\cal R}^2\pi_{\ast}\bigl({\cal E}_{e_{II}}\bigr)$, 
because by relative Serre vanishing theorem 
${\cal R}^i\pi_{\ast}\bigl({\cal O}_{\cal X}(nD)\otimes {\cal E}_{e_{II}}\bigr)=0$
 for large enough $n$.

 By our simplifying assumption \ref{assum; 2} on $e_{II}$ (on page 
\pageref{assum; 2}), 
${\cal R}^2\pi_{\ast}\bigl({\cal E}_{e_{II}}\bigr)=0$. So 
 ${\cal R}^1\pi_{\ast}\bigl({\cal O}_{nD}(nD)\otimes {\cal E}_{e_{II}}\bigr)=0$
 and ${\cal R}^0\pi_{\ast}\bigl({\cal O}_{nD}(nD)\otimes {\cal E}_{e_{II}}\bigr)$
 is locally free.  Then we may take ${\cal V}_{e_{II}}$ 
and ${\cal W}_{e_{II}}$ to be
 the locally free sheaves 
${\cal R}^0\pi_{\ast}\bigl({\cal O}_{\cal X}(nD)\otimes {\cal E}_{e_{II}}
\bigr)$ and
 ${\cal R}^0\pi_{\ast}\bigl({\cal O}_{nD}(nD)\otimes {\cal E}_{e_{II}}\bigr)$,
 respectively.

\medskip

 Set ${\bf V}_{e_{II}}$, ${\bf W}_{e_{II}}$ to be the vector bundles
associated with the locally free sheaves 
${\cal V}_{e_{II}}={\cal R}^0\pi_{\ast}\bigl({\cal O}_{\cal X}(nD)\otimes
 {\cal E}_{e_{II}}\bigr)$ and ${\cal W}_{e_{II}}=
{\cal R}^0\pi_{\ast}\bigl({\cal O}_{nD}(nD)
\otimes {\cal E}_{e_{II}}\bigr)$, respectively.
 These bundles depend on the choices of $D$ and $n$ and are not canonical\footnote{
We drop the notational dependence of 
 ${\bf V}_{e_{II}}$ and ${\bf W}_{e_{II}}$ 
on $D$ and $n$ to simplify our symbols.}. 

Then $\Phi_{{\bf V}_{e_{II}}{\bf W}_{e_{II}}}:{\bf V}_{e_{II}}\mapsto 
{\bf W}_{e_{II}}$ defines an algebraic family Kuranishi model of $e_{II}$,
 as was described in section \ref{subsection; type2K}. 
Under the assumption that $e_{II}-c_1({\bf K}_{{\cal X}/B})$ is nef, 
by Riemann-Roch theorem 
$rank_{\bf C}({\bf V}_{e_{II}}-{\bf W}_{e_{II}})=1-q+p_g+{e_{II}^2-
c_1({\bf K}_{{\cal X}/B})\cdot e_{II}\over 2}$.

\begin{rem}\label{rem; geometric}
It is useful to comprehend the
 geometric meaning of the above Kuranishi model. By adjoining
 a very ample $nD$, the family moduli space of curves  ${\cal M}_{e_{II}}$ is
 naturally
 embedded into the family moduli space of a better behaved class 
$e_{II}+nD$, which \footnote{Thanks
 to the sufficiently very 
ampleness of $D$ and the large number $n$.} has the nice structure of
 a projective space bundle over ${\cal T}_B{\cal X}$. Then the sub-locus 
 ${\cal M}_{e_{II}}$ is
 characterized by the cross section of the obstruction bundle induced by
 $\Phi_{{\bf V}_{e_{II}}{\bf W}_{e_{II}}}$, which requries the ray of
 non-zero 
sections in $H^0({\cal X}_b, {\cal E}_{e_{II}}\otimes {\cal O}_{\cal X}(nD)|_b)$
 to vanish along the effective divisor $nD$ to recover a curve representing
 $e_{II}$.
\end{rem}

 In the above argument, we have not made use of the exceptionality property, 
i.e. definition \ref{defin; exception}, on
 $e_{II}$. So we may replace $e_{II}$ by any $\underline{C}$ or 
$\underline{C}-e_{II}$ which satisfies the nef condition
 on $\underline{C}-e_{II}-c_1({\bf K}_{{\cal X}/B})$. As before,
 we still assume $febd(\underline{C}, {\cal X}/B)=febd(\underline{C}-e_{II}, 
{\cal X}/B)=p_g$ to simplify our discussion.
 The above argument still goes through without modification. 
 We can choose the same effective ample $D$ and a uniformly large $n$ such that the
sheaf morphisms
${\cal R}^0\pi_{\ast}\bigl({\cal O}_{\cal X}(nD)\otimes
 {\cal E}_{\underline{C}}\bigr)\mapsto 
{\cal R}^0\pi_{\ast}\bigl({\cal O}_{nD}(nD)
\otimes {\cal E}_{\underline{C}}\bigr)$
and ${\cal R}^0\pi_{\ast}\bigl({\cal O}_{\cal X}(nD)\otimes
 {\cal E}_{\underline{C}-e_{II}}\bigr)\mapsto 
{\cal R}^0\pi_{\ast}\bigl({\cal O}_{nD}(nD)
\otimes {\cal E}_{\underline{C}-e_{II}}\bigr)$ define the
algebraic family Kuranishi models of $\underline{C}$ and $\underline{C}-e_{II}$,
 respectively.

We denote the corresponding Kuranishi-model vector bundles by
 ${\bf V}_{\underline{C}}$, ${\bf W}_{\underline{C}}$ and
 ${\bf V}_{\underline{C}-e_{II}}$, ${\bf W}_{\underline{C}-e_{II}}$, respectively.
 In the following, we will fix a pair of $D$ and $n$ and discuss the
 switching of the family Kuranishi models between $\underline{C}$ and 
$\underline{C}-e_{II}$.
 
\begin{rem}\label{rem; similar}
 We notice that if we formally replace $nD+\underline{C}$ by $C$, $nD$ by 
${\bf M}(E)E$ and $\underline{C}$ by $C-{\bf M}(E)E$, 
the above algebraic family Kuranishi model of $\underline{C}$
 corresponds formally 
to the canonical algebraic family Kuranishi model of $C-{\bf M}(E)E$, introduced
 in [Liu3] and used heavily in [Liu6],

$${\cal R}^0\pi_{\ast}\bigl({\cal O}_{M_{n+1}}\otimes {\cal E}_C\bigr)\mapsto
 {\cal R}^0\pi_{\ast}\bigl({\cal O}_{{\bf M}(E)E}\otimes {\cal E}_C\bigr).$$
\end{rem}

 This analogue provides us an easy way to memorize and link their
 family Kuranishi models.

\bigskip

\subsection{The Switching of Family Kuranishi Models Involving type $II$
 Exceptional Classes}\label{subsection; type II}

\bigskip

 In the following, we compare the Kuranishi datum of $\underline{C}$ and
 $\underline{C}-e_{II}$ (using only one $e_{II}$).
 Consider the pull-back of the Kuranishi models
 datum of $\underline{C}$ and $\underline{C}-e_{II}$ 
from $B$ to ${\cal M}_{e_{II}}$ by the
 natural projection ${\cal M}_{e_{II}}\mapsto B$.

 Let ${\bf e}_{II}$, with the following commutative diagram,

\[
\begin{array}{ccc}
{\bf e}_{II} & \hookrightarrow & {\cal X}\times_B{\cal M}_{e_{II}}\\
\Big\downarrow & & \Big\downarrow \\
{\cal M}_{e_{II}} & = & {\cal M}_{e_{II}} 
\end{array}
\]

denote the universal type $II$ 
curve over ${\cal M}_{e_{II}}$.

\begin{prop}\label{prop; compare}
 Consider the pull-backs of the family Kuranishi models of 
$\underline{C}$ and $\underline{C}-e_{II}$ to ${\cal M}_{e_{II}}$.
Between the pull-backs to ${\cal M}_{e_{II}}$ of the
algebraic family Kuranishi models of $\underline{C}$ and
 $\underline{C}-e_{II}$ constructed following the recipe 
in subsection \ref{subsection; type2K},
 there is a commutative diagram of
 ``columns of short exact sequence'' of locally free sheaves,

\[
\begin{array}{ccc}
0 & & 0\\
\Big\downarrow & & \Big\downarrow\\
{\cal R}^0\pi_{\ast}\bigl({\cal O}_{\cal X}(nD)\otimes 
{\cal E}_{\underline{C}-e_{II}}\bigr)&\stackrel{\Phi_{{\cal V}_{\underline{C}
-e_{II}}{\cal W}_{\underline{C}-e_{II}}}}{\longrightarrow}
 & {\cal R}^0\pi_{\ast}\bigl({\cal O}_{nD}(nD)\otimes 
{\cal E}_{\underline{C}-e_{II}}\bigr)\\
\Big\downarrow & & \Big\downarrow\\
{\cal R}^0\pi_{\ast}\bigl({\cal O}_{\cal X}(nD)\otimes 
{\cal E}_{\underline{C}}\bigr)\otimes {\cal H}_{II} & 
\stackrel{\Phi_{{\cal V}_{\underline{C}}{\cal 
W}_{\underline{C}}}\otimes 
id_{{\cal H}_{II}}}{\longrightarrow}& {\cal R}^0\pi_{\ast}\bigl(
{\cal O}_{nD}(nD)\otimes 
{\cal E}_{\underline{C}}\bigr) \otimes {\cal H}_{II}\\
\Big\downarrow & & \Big\downarrow\\
{\cal R}^0\pi_{\ast}\bigl({\cal O}_{{\bf e}_{II}}(nD)
\otimes {\cal E}_{\underline{C}}\bigr) \otimes {\cal H}_{II} & 
\stackrel{}{\longrightarrow}& {\cal R}^0\pi_{\ast}\bigl({\cal O}_{{\bf e}_{II}\cap
 nD}(nD)\otimes {\cal E}_{\underline{C}}\bigr)\otimes {\cal H}_{II} \\
\Big\downarrow & & \Big\downarrow\\
0 & & 0\\
\end{array}
\]
\end{prop}

 Here the hyperplane bundle ${\cal H}_{II}\mapsto {\cal M}_{e_{II}}$ in 
the above diagram is 
induced from the embedding ${\cal M}_{e_{II}}\hookrightarrow 
{\bf P}({\bf V}_{e_{II}})$.

\medskip

\noindent Sketch of the Proof: 
 Both of the columns are the derived long exact sequences
 from some twisted versions of the following short exact sequence, pull-back to
 ${\cal X}\times_B{\cal M}_{e_{II}}\mapsto {\cal M}_{e_{II}}$,

$
0\mapsto {\cal O}_{\cal X}(-{\bf e}_{II})\mapsto {\cal O}_{\cal X}\mapsto 
{\cal O}_{{\bf e}_{II}}\mapsto 0$

 and its restriction to the non-reduced $nD$. These derived sequences
 are short exact as $n$ has been
 chosen to be large enough to guarantee the vanishing of
 ${\cal R}^i\pi_{\ast}\bigl({\cal O}_{\cal X}(nD)\otimes {\cal E}\bigr)$,
 $i>0$ for either ${\cal E}={\cal E}_{\underline{C}}$ or 
${\cal E}_{\underline{C}-e_{II}}$.
 The commutativity of the sequences follow from the commutativity of
 the tensoring operations
 (by tensoring the defining sections of ${\bf e}_{II}\subset 
{\cal X}\times_B{\cal M}_{e_{II}}$) and the restriction
 to $nD$. The map of the last row is induced by 
the derived exact sequence of
 $0\mapsto {\cal O}_{{\bf e}_{II}}\otimes {\cal E}_{\underline{C}}\mapsto 
{\cal O}_{{\bf e}_{II}}(nD)\otimes {\cal E}_{\underline{C}}\mapsto
 {\cal O}_{{\bf e}_{II}\cap nD}(nD)\otimes {\cal E}_{\underline{C}}\mapsto 0$.

 If ${\bf e}_{II}\cap nD\mapsto {\cal M}_{e_{II}}$ has been a finite morphism,
 the sheaf ${\cal R}^0\pi_{\ast}\bigl({\cal O}_{{\bf e}_{II}\cap
 nD}(nD)\otimes {\cal E}_{\underline{C}}\bigr)$ is automatically locally free 
with its rank equal to the relative length of ${\bf e}_{II}\cap nD\mapsto
 {\cal M}_{e_{II}}$. Without the finiteness assumption of the morphism
${\bf e}_{II}\cap nD\mapsto {\cal M}_{e_{II}}$, we check the locally freeness
 of ${\cal R}^0\pi_{\ast}\bigl({\cal O}_{{\bf e}_{II}\cap
 nD}(nD)\otimes {\cal E}_{\underline{C}}\bigr)$ by proving that the
second column in the above commutative diagram remains left exact after 
tensoring with $k(y)$, for all the closed points $y\in {\cal M}_{e_{II}}$.
$\Box$

\bigskip

  The exact sequences in proposition \ref{prop; compare}
 imply that ${\bf V}_{\underline{C}-e_{II}}$ is
 a sub-bundle of ${\bf V}_{\underline{C}}\otimes {\bf H}_{II}$.
 Denote the quotient bundles
 associated with 
${\cal R}^0\pi_{\ast}\bigl({\cal O}_{e_{II}}(nD)\otimes {\cal E}_{\underline{C}}
\bigr)\otimes {\cal H}_{II}$ 
and ${\cal R}^0\pi_{\ast}\bigl({\cal O}_{{\bf e}_{II}\cap
 nD}(nD)\otimes {\cal E}_{\underline{C}}\bigr)\otimes {\cal H}_{II}$
by ${\bf V}'$ and ${\bf W}'$, respectively.
 Then ${\bf P}_B({\bf V}_{\underline{C}-e_{II}})$ can be
 viewed as the smooth sub-scheme of 
${\bf P}({\bf V}_{\underline{C}}\otimes {\bf H}_{II})$
$\cong {\bf P}({\bf V}_{\underline{C}})$
 defined by the zero locus of the section of ${\bf H}\otimes 
\pi_{{\bf P}_B({\bf V}_{\underline{C}})}^{\ast}{\bf V}'$ induced by
 ${\bf V}_{\underline{C}}\otimes {\bf H}_{II}\mapsto {\bf V}'\mapsto 0$.

 This implies that we may replace the original ambient space of the
 family moduli space ${\cal M}_{\underline{C}-e_{II}}$ 
of $\underline{C}-e_{II}$ from
 ${\bf P}_B({\bf V}_{\underline{C}-e_{II}})$ by
 ${\bf P}_B({\bf V}_{\underline{C}})$ and replace the original
 obstruction bundle ${\bf H}\otimes 
\pi_{{\bf P}_B({\bf V}_{\underline{C}}\otimes 
{\bf H}_{II})}^{\ast}{\bf W}_{\underline{C}-e_{II}}$
 by an extended obstruction bundle equivalent to ${\bf H}\otimes 
\pi_{{\bf P}_B({\bf V}_{\underline{C}\otimes 
{\bf H}_{II}})}^{\ast}({\bf W}_{\underline{C}-
e_{II}}\oplus {\bf V}')$ in the appropriated $K$ group.

\bigskip

 In the following we discuss how the desired 
``extended obstruction bundle'' can be constructed
from the standard bundle extension construction.

\medskip

 Consider a bundle extension of ${\bf V}'$ by
 ${\bf W}_{\underline{C}-e_{II}}$. All such bundle extensions are classified
 by the group $Ext^1({\bf V}', {\bf W}_{\underline{C}-e_{II}})$.
 By applying the left exact functor $HOM({\bf V}', \bullet)$ to
 the short exact sequence $0\mapsto {\bf W}_{\underline{C}-e_{II}}\mapsto
 {\bf W}_{\underline{C}}\otimes {\bf H}_{II}\mapsto {\bf W}'
\mapsto 0$, we get the
 following portion of derived long exact sequence,

$$ HOM({\bf V}', {\bf W}_{\underline{C}}\otimes {\bf H}_{II})\mapsto
 HOM({\bf V}', {\bf W}')\stackrel{\delta}{\longrightarrow}
Ext^1({\bf V}', {\bf W}_{\underline{C}-e_{II}})\cdots.$$

 Apparently the bundle map ${\bf V}'\mapsto {\bf W}'$ induced from 
 the sheaf commutative diagram in proposition \ref{prop; compare} gives 
an element in $HOM({\bf V}', {\bf W}')$. Its image in 
$Ext^1({\bf V}', {\bf W}_{\underline{C}-e_{II}})$ under the
 connecting homomorphism determines a bundle
extension and therefore defines the new extended 
obstruction bundle ${\bf W}_{new}$. 
 And we have the following defining short exact sequence 

$$0\mapsto {\bf W}_{\underline{C}-e_{II}}\mapsto {\bf W}_{new}\mapsto
 {\bf V}'\mapsto 0.$$

 To show that there is a canonically induced 
bundle map ${\bf W}_{new}\mapsto {\bf W}_{\underline{C}}\otimes {\bf H}_{II}$, 
we take
 $HOM({\bf W}_{new}, \bullet)$ on the short exact sequence 
$0\mapsto {\bf W}_{\underline{C}-e_{II}}\mapsto 
{\bf W}_{\underline{C}}\otimes {\bf H}_{II}\mapsto
 {\bf W}'\mapsto 0$ and get the following portion of derived long
exact sequence,

$$ \cdots\mapsto HOM({\bf W}_{new}, {\bf W}_{\underline{C}}\otimes {\bf H}_{II})
\mapsto HOM({\bf W}_{new}, {\bf W}')\mapsto Ext^1({\bf W}_{new}, 
{\bf W}_{\underline{C}-e_{II}})\mapsto \cdots.$$

 The composition ${\bf W}_{new}\mapsto {\bf V}'\mapsto {\bf W}'$ induces
an element in $HOM({\bf W}_{new}, {\bf W}')$. To show that
 it is the image of an element in $HOM({\bf W}_{new}, 
{\bf W}_{\underline{C}}\otimes {\bf H}_{II})$,
 it suffices to show that its image into 
$Ext^1({\bf W}_{new}, {\bf W}_{\underline{C}-e_{II}})$ vanishes.

\medskip

On the other hand, the derived long exact sequence of the contravariant 
functor 
$HOM(\bullet, {\bf W}_{\underline{C}-e_{II}})$ upon the defining 
short exact sequence of ${\bf W}_{new}$,
 $0\mapsto {\bf W}_{\underline{C}-e_{II}}\mapsto {\bf W}_{new}\mapsto 
{\bf V}'\mapsto 0$ implies that

$$HOM({\bf W}_{\underline{C}-e_{II}}, {\bf W}_{\underline{C}-e_{II}})
\mapsto Ext^1({\bf V}', {\bf W}_{\underline{C}-e_{II}})
\mapsto Ext^1({\bf W}_{new}, {\bf W}_{\underline{C}-e_{II}})\cdots.$$

 The extension class 
$\in Ext^1({\bf V}', {\bf W}_{\underline{C}-e_{II}})$ defining
 ${\bf W}_{new}$ is the
 image of $id_{{\bf W}_{\underline{C}-e_{II}}}
\in HOM({\bf W}_{\underline{C}-e_{II}}, {\bf W}_{\underline{C}-e_{II}})$.
 Therefore by 
exactness of the above derived sequence its image in 
$Ext^1({\bf W}_{new}, {\bf W}_{\underline{C}-e_{II}})$ must vanish.

\bigskip

 The bundle morphism 
 $\in HOM({\bf W}_{new}, {\bf W}_{\underline{C}}\otimes {\bf H}_{II})$ 
restricts to the the bundle injection 
${\bf W}_{\underline{C}-e_{II}}\mapsto {\bf W}_{\underline{C}}\otimes 
{\bf H}_{II}$ on the
sub-bundle ${\bf W}_{\underline{C}-e_{II}}\subset {\bf W}_{new}$.

 To show that ${\bf V}_{\underline{C}}\otimes {\bf H}_{II}\mapsto {\bf V}'$ 
can be lifted \footnote{Notice that the lifting may not be unique!}
to some 
${\bf V}_{\underline{C}}\otimes {\bf H}_{II}\mapsto {\bf W}_{new}$, it suffices
 to show that its image $\iota$ in 
$Ext^1({\bf V}_{\underline{C}}\otimes {\bf H}_{II}, 
{\bf W}_{\underline{C}-e_{II}})$ vanishes. This is 
ensured by the following derived long exact sequence,

$$\cdots HOM({\bf V}_{\underline{C}}\otimes {\bf H}_{II}, 
{\bf W}_{\underline{C}}\otimes {\bf H}_{II})
\mapsto HOM({\bf V}_{\underline{C}}\otimes {\bf H}_{II}, {\bf W}')\mapsto 
Ext^1({\bf V}_{\underline{C}}\otimes {\bf H}_{II}, 
{\bf W}_{\underline{C}-e_{II}})\cdots,$$

  the fact that $\iota$ is the composite image of 
 the element ${\bf V}_{\underline{C}}\otimes {\bf H}_{II}
\mapsto {\bf W}_{\underline{C}}\otimes {\bf H}_{II}$ in
 $HOM({\bf V}_{\underline{C}}\otimes {\bf H}_{II}, 
{\bf W}_{\underline{C}}\otimes {\bf H}_{II})$, and the
acyclicity of the above long exact sequence.

\medskip

The following lemma fixes the unique lifting of ${\bf V}_{\underline{C}}\otimes
 {\bf H}_{II}\mapsto {\bf V}'$.

\begin{lemm}\label{lemm; compatible}
Among the possible liftings of ${\bf V}_{\underline{C}}\otimes {\bf H}_{II}
\mapsto {\bf V}'$,
there is a unique lifting ${\bf V}_{\underline{C}}\otimes {\bf H}_{II}\mapsto 
{\bf W}_{new}$ which makes the following diagram commutative,

\[
\begin{array}{ccc}
 {\bf V}_{\underline{C}}\otimes {\bf H}_{II} & \longrightarrow & {\bf W}_{new}\\
 \Big\downarrow & & \Big\downarrow \\
 {\bf V}_{\underline{C}}\otimes {\bf H}_{II} 
& \longrightarrow &  {\bf W}_{\underline{C}}\otimes {\bf H}_{II}
\end{array}
\]

\end{lemm}

The above commutative diagram ensures the compatibility between the 
(extended) algebraic family
Kuranishi models of $\underline{C}-e_{II}$ and of $\underline{C}$ above
 ${\cal M}_{e_{II}}$.

\medskip

\noindent Proof of Lemma \ref{lemm; compatible}:
Start with an arbitrary lifting $\kappa:{\bf V}_{\underline{C}}\otimes 
{\bf H}_{II}\mapsto 
{\bf W}_{new}$. We have
 the following commutative diagram among $HOM$ groups,

\[
\hskip -.4in
\begin{array}{ccccc}
HOM({\bf V}_{\underline{C}}\otimes {\bf H}_{II}, 
{\bf W}_{\underline{C}-e_{II}}) & \mapsto 
& HOM({\bf V}_{\underline{C}}\otimes {\bf H}_{II}, {\bf W}_{new} )& \mapsto & 
HOM({\bf V}_{\underline{C}}\otimes {\bf H}_{II}, {\bf V}')\\
\Big\downarrow \vcenter{%
\rlap{$\scriptstyle{\mathrm{=}}\,$}}
& & \Big\downarrow  & & \Big\downarrow \\
HOM({\bf V}_{\underline{C}}\otimes {\bf H}_{II}, 
{\bf W}_{\underline{C}-e_{II}})& \mapsto 
& HOM({\bf V}_{\underline{C}}\otimes {\bf H}_{II}, 
{\bf W}_{\underline{C}}\otimes {\bf H}_{II})& \mapsto &
HOM({\bf V}_{\underline{C}}\otimes {\bf H}_{II}, {\bf W}')
\end{array}
\]

 We consider the difference between
 ${\bf V}_{\underline{C}}\otimes {\bf H}_{II}\mapsto
 {\bf W}_{\underline{C}}\otimes {\bf H}_{II}$ and the image of $\kappa$ in
 $HOM({\bf V}_{\underline{C}}\otimes {\bf H}_{II}, 
{\bf W}_{\underline{C}}\otimes {\bf H}_{II})$ and denote it
by $\zeta$. The image of $\zeta$ in 
$HOM({\bf V}_{\underline{C}}\otimes {\bf H}_{II}, {\bf W}')$
 vanishes because both the image of $\kappa$ and 
${\bf V}_{\underline{C}}\otimes {\bf H}_{II}\mapsto
 {\bf W}_{\underline{C}}\otimes {\bf H}_{II}$ have the same image
 ${\bf V}_{\underline{C}}\otimes {\bf H}_{II}\mapsto {\bf W}'$ in  
$HOM({\bf V}_{\underline{C}}\otimes {\bf H}_{II}, {\bf W}')$. By exactness of the
 second row above,  there exists an element $\rho\in 
HOM({\bf V}_{\underline{C}}\otimes {\bf H}_{II}, 
{\bf W}_{\underline{C}-e_{II}})$ which maps
 onto $\zeta$. By the commutativity of the diagram, denote 
the image of the element 
$\rho$ in $HOM({\bf V}_{\underline{C}}\otimes {\bf H}_{II}, 
{\bf W}_{new})$ by $\rho'$.
Then we replace $\kappa$ by $\kappa+\rho'$ and it is rather easy to see that
 the image of $\kappa+\rho'$ in 
$HOM({\bf V}_{\underline{C}}\otimes {\bf H}_{II}, 
{\bf W}_{\underline{C}}\otimes {\bf H}_{II})$ is
 ${\bf V}_{\underline{C}}\otimes {\bf H}_{II}\mapsto 
{\bf W}_{\underline{C}}\otimes {\bf H}_{II}$.
 The commutativity of the original diagram follows.
$\Box$

\bigskip

 It is vital to understand the bundle map ${\bf W}_{new}\mapsto 
{\bf W}_{\underline{C}}\otimes {\bf H}_{II}$. If the bundle map is injective
 over ${\cal M}_{e_{II}}$, then
 the restriction of the 
family moduli space of $\underline{C}-e_{II}$ is isomorphic
 to the restriction of family moduli space of $\underline{C}$ 
above ${\cal M}_{e_{II}}$.
 Namely, we have the isomorphism
 ${\cal M}_{\underline{C}-e_{II}}\times_B{\cal M}_{e_{II}}\cong
 {\cal M}_{\underline{C}}\times_B{\cal M}_{e_{II}}$. This will be formulated
 as a {\bf special condition} on page \pageref{special}.

 On the other hand, the possible failure of the injectivity of 
${\bf W}_{new}\mapsto {\bf W}_{\underline{C}}\otimes {\bf H}_{II}$ over
 ${\cal M}_{e_{II}}$ may result in
 the discrepancy of identifying the fiber products 
${\cal M}_{\underline{C}-e_{II}}\times_B{\cal M}_{e_{II}}$ with 
 ${\cal M}_{\underline{C}}\times_B{\cal M}_{e_{II}}$.

In the following, we identify the kernel cone of
 ${\bf W}_{new}\mapsto {\bf W}_{\underline{C}}\otimes {\bf H}_{II}$.

\begin{prop}\label{prop; thesame}
Given the unique bundle map lifting ${\bf V}_{\underline{C}}\otimes {\bf H}_{II}
\mapsto {\bf W}_{new}$ of ${\bf V}_{\underline{C}}\otimes {\bf H}_{II}\mapsto
 {\bf V}'$
fixed by lemma \ref{lemm; compatible}, 
the following commutative diagram of vector bundles

\[
\begin{array}{ccc}
{\bf W}_{new} & \longrightarrow & {\bf W}_{\underline{C}}\otimes {\bf H}_{II}\\
\Big\downarrow & & \Big\downarrow\\
{\bf V}' & \longrightarrow & {\bf W}'
\end{array}
\]

induces an isomorphism between the kernel cones of the
horizontal bundle morphisms.
\end{prop}

\noindent Proof of proposition \ref{prop; thesame}: 
The vertical arrows are known to be bundle
surjections by our construction. 
The kernels of both the vertical bundle maps of the
 columns are isomorphic to
 ${\bf W}_{\underline{C}-e_{II}}$. Let ${\bf C}_{new}$ and $ {\bf C}'$ be
the kernel sub-cones of ${\bf W}_{new}\mapsto 
{\bf W}_{\underline{C}}\otimes {\bf H}_{II}$
 and of ${\bf V}'\mapsto {\bf W}'$, respectively. We prove that
 ${\bf C}_{new}\mapsto {\bf C}'$ induced from ${\bf W}_{new}\mapsto {\bf V}'$
  is an isomorphism of abelian 
cones.  , i.e., 
it suffices to show that for all the closed $b\in B$ the
 fibers of the cones are isomorphic under ${\bf W}_{new}|_b\mapsto
{\bf V}'|_b$. Because the restriction of 
 ${\bf W}_{new}\mapsto {\bf W}_{\underline{C}}\otimes {\bf H}_{II}$ to 
the sub-bundle ${\bf W}_{\underline{C}-e_{II}}$
 has been injective, so ${\bf C}_{new}
\cap {\bf W}_{\underline{C}-e_{II}}$ is embedded
 as the
 zero section sub-cone of ${\bf W}_{new}$.
By the above commutative diagram it is
clear that ${\bf C}_{new}\mapsto {\bf C}'$.
 On the other hand for all vectors $t\in {\bf C}'|_b$, an arbitrary
lifting of $t$, $\tilde{t}\in {\bf W}_{new}|_b$, 
may or may not lie in the kernel space above $b$, ${\bf C}_{new}|_b$.
 But the image of $t\in {\bf C}'|_b$ into ${\bf W}'|_b$ is zero. So the image of
 $\tilde{t}$ in ${\bf W}_{\underline{C}}\otimes {\bf H}_{II}$ 
has to map trivially into ${\bf W}'|_b$.
So one may find a unique
 element $w(t)$ in ${\bf W}_{\underline{C}-e_{II}}|_b$ which maps
 onto the image of $\tilde{t}$ in
 ${\bf W}_{\underline{C}}\otimes {\bf H}_{II}|_b$. 
Because ${\bf W}_{\underline{C}-e_{II}}\subset {\bf W}_{new}$, 
 $w(t)$ can be viewed as an element
 in ${\bf W}_{new}|_b$ as well and then $t-w(t)$ will be mapped trivially into 
 ${\bf W}_{\underline{C}}\otimes {\bf H}_{II}|_b$. So $t-w(t)\in {\bf C}_{new}|_b$.
We have shown that every element 
 $t\in {\bf C}'$ can be lifted uniquely into an element in ${\bf C}_{new}|_b$
 and we establish their
bijections for all the closed points $b\in B$.  
$\Box$

\begin{rem}\label{rem; many}
 In the above discussion we have hardly used any special property of $e_{II}$. 
Suppose that there are $p$ distinct type $II$ exceptional classes 
$e_{II; 1}$, $e_{II; 2}$, $\cdots$, $e_{II; p}$ over ${\cal X}\mapsto B$,
 with $e_{II; i}\cdot e_{II; j}\geq 0$ for $i\not=j$.
If we replace ${\cal M}_{e_{II}}$ by the co-existence locus \footnote{
It will be defined and discussed in more details in subsection 
\ref{subsection; stable}.}
 of the
 type II classes, $\times_B^{i\leq p} {\cal M}_{e_{II; i}}$,
 and replace the single universal curve ${\bf e}_{II}$ by the total sum
 $\sum_{i\leq p} {\bf e}_{II; i}$ of all the universal curves, 
the above discussion can be generalized
easily to the cases involving more than one type $II$ class.
\end{rem}

\begin{rem}\label{rem; analogue}
The above discussion and remark \ref{rem; many} 
imply that ${\bf W}_{new}\mapsto {\bf W}_{\underline{C}}\otimes {\bf H}_{II}$
 has been the analogue of ${\bf W}_{canon}^{\circ}\mapsto {\bf W}_{canon}$ in 
comparing the family Kuranishi models of $C-{\bf M}(E)E-\sum e_{k_i}$ and
 $C-{\bf M}(E)E$ involving type $I$ exceptional classes.
The main difference is that in the type $II$ case they are not canonical.
\end{rem}

 Recall that in [Liu5], [Liu6], we have discussed the relationship of
 the canonical algebraic family Kuranishi models of 
$C-{\bf M}(E)E$ and $C-{\bf M}(E)E-\sum_{i\leq p}
 e_{k_i}$ over $Y(\Gamma)\times T(M)$, 
where $e_{k_i}\cdot (C-{\bf M}(E)E)<0$ for $1\leq i\leq p$.

Let $\Phi_{{\bf V}_{canon}^{\circ}{\bf W}_{canon}^{\circ}}:
{\bf V}_{canon}^{\circ}\mapsto {\bf W}_{canon}^{\circ}$ and 
$\Phi_{{\bf V}_{canon}{\bf W}_{canon}}:
{\bf V}_{canon}\mapsto {\bf W}_{canon}$ be the 
canonical algebraic family Kuranishi models of 
$C-{\bf M}(E)E$ and $C-{\bf M}(E)E-\sum e_{k_i}$, respectively.
Please consult lemma 6 of [Liu5] for their definitions. 

 Then we have ${\bf V}_{canon}^{\circ}={\bf V}_{canon}$ and
 the bundle map over 
$Y(\Gamma)\times T(M)$, ${\bf W}_{canon}^{\circ}|_{Y(\Gamma)\times T(M)}
\mapsto {\bf W}_{canon}|_{Y(\Gamma)\times T(M)}$, parallel to the
 sheaf sequence below.

$$\hskip -1.2in {\cal R}^0\pi_{\ast}\bigl({\cal O}_{\sum \Xi_{k_i}}\otimes
 {\cal E}_{C-{\bf M}(E)E}\bigr)\mapsto
{\cal R}^0\pi_{\ast}\bigl({\cal O}_{{\bf M}(E)E+\sum \Xi_{k_i}}\otimes {\cal E}_C
\bigr)\mapsto {\cal R}^0\pi_{\ast}\bigl(
{\cal O}_{{\bf M}(E)E}\otimes {\cal E}_C\bigr)\mapsto
{\cal R}^1\pi_{\ast}\bigl({\cal O}_{\sum \Xi_{k_i}}\otimes
 {\cal E}_{C-{\bf M}(E)E}\bigr),$$

 for the sum of the type $I$ universal curves $\Xi_{k_i}\mapsto Y(\Gamma)$.

 It is not hard to establish the following correspondence,

\medskip

\noindent $\bullet$ ${\bf V}_{canon}^{\circ}={\bf V}_{canon}$ corresponds to
 ${\bf V}_{\underline{C}}\otimes {\bf H}_{II}=
{\bf V}_{\underline{C}}\otimes {\bf H}_{II}$.

\medskip

\noindent $\bullet$ ${\bf W}_{canon}^{\circ}\mapsto {\bf W}_{canon}$ corresponds 
to ${\bf W}_{new}\mapsto {\bf W}_{\underline{C}}\otimes {\bf H}_{e_{II}}$.

\medskip

\begin{rem}\label{rem; independent}
 For the purpose of defining the virtual fundamental classes
 of ${\cal M}_{\underline{C}}$, it is easy to
see that the twisting operation from 
${\bf V}_{\underline{C}}\mapsto {\bf W}_{\underline{C}}$ to
${\bf V}_{\underline{C}}\otimes {\bf H}_{II}
\mapsto {\bf W}_{\underline{C}}\otimes {\bf H}_{II}$ does not affect the
 virtual fundamental class $[{\cal M}_{\underline{C}}]_{vir}$
\footnote{The line bundle ${\bf H}_{II}$ does not show up in the theory of the
type $I$ exceptional classes. Since for type $I$ classes we
 use the canonical algebraic family Kuranishi models of $e_{k_i}$ and
${\bf H}_{II}$ is
 reduced to trivial line bundle over 
${\cal M}_{e_{k_i}}\cong Y(\Gamma_{e_{k_i}})$.} of ${\cal M}_{\underline{C}}$.
\end{rem}
 
\medskip

\subsection{The Vector Bundle ${\bf W}'$ and the Type $II$ Class $e_{II}$
}\label{subsection; bundle}

\bigskip

 In the previous subsection, section \ref{subsection; type II}, 
we have defined ${\bf W}'$ to be the vector bundle associated with
${\cal R}^0\pi_{\ast}\bigl({\cal O}_{{\bf e}_{II}\cap nD}(nD)\otimes
 {\cal E}_{\underline{C}}\bigr)\otimes {\cal H}_{II}$ and it is a 
quotient bundle of ${\bf W}_{\underline{C}}\otimes {\bf H}_{II}$. 
In this subsection, we
 would like to bridge ${\bf W}'$ with the virtual fundamental class of the
 type $II$ class $e_{II}$. The discussion can be extended to
 more than one $e_{II; i}$ in a parallel manner, once we substitute
 $e_{II}$ and ${\bf e}_{II}$ by $e_{II; i}$ and $\sum {\bf e}_{II; i}$, 
respectively. The general version will be address in subsection 
\ref{subsection; virtual} later.
 
Firstly we prove a simple lemma,

\begin{lemm}\label{lemm; factor}
 Let ${\bf W}_{II}'\mapsto {\cal M}_{e_{II}}$ denote the vector bundle
 associated with the locally free sheaf ${\cal R}^0\pi_{\ast}\bigl(
{\cal O}_{{\bf e}_{II}\cap nD}({\bf e}_{II}+nD)\bigr)\otimes {\cal H}_{II}$. 
Then we have the equality on the ranks 
$rank_{\bf C}{\bf W}_{II}'=rank_{\bf C}{\bf W}'$,
and the total Chern class of ${\bf W}'$, $c_{total}({\bf W}')$, can be 
identified with $c_{total}({\bf W}_{II}')+\eta$, where $\eta$ is 
a polynomial of the ``${\bf e}_{II}\mapsto {\cal M}_{e_{II}}$'' push-forward of 
monomials the variables
 $c_1({\cal E}_{\underline{C}-e_{II}})$, $e_{II}$ and $nD$.
\end{lemm}

\noindent Proof: To determine the ranks of ${\bf W}_{II}'$ and ${\bf W}'$,
 it suffices to calculate them at a closed point $b\in {\cal M}_{e_{II}}$.

Because $nD$ is very ample in ${\cal X}$, we can replace $nD$ by a
 linearly equivalent very ample divisor such that it intersects with the
 curve ${\bf e}_{II}|_b$ at a finite number of points. It is easy to see
 by base-change theorem that the ranks of both ${\bf W}_{II}'$ and ${\bf W}'$
 are equal to $\int_{\cal X}{\bf e}_{II}|_b\cap nD\in {\bf N}$. Thus they
 must be equal.

 Because the higher derived images vanish, both
 ${\cal R}^0\pi_{\ast}({\cal O}_{{\bf e}_{II}\cap nD}({\bf e}_{II}+nD))=\pi_{\ast}
({\cal O}_{{\bf e}_{II}\cap nD}({\bf e}_{II}+nD))$ and
 ${\cal R}^0\pi_{\ast}({\cal O}_{{\bf e}_{II}\cap nD}(\underline{C}+nD))=\pi_{\ast}
({\cal O}_{{\bf e}_{II}\cap nD}(\underline{C}+nD))$, and their
 Chern characters can be computed by Grothendieck-Riemann-Roch formula (See
 chapter 15. and chapter 18. of [F]).

 Since ${\cal O}_{{\bf e}_{II}\cap nD}(\underline{C}+nD)$ can be
 constructed from ${\cal O}_{{\bf e}_{II}\cap nD}({\bf e}_{II}+nD)$ by
twisting ${\cal E}_{\underline{C}-e_{II}}$, the conclusion follows from
 Grothendieck-Riemann-Roch formula and the fact that total Chern class
 can be expressed in terms of the Chern character algebraically.
$\Box$

\medskip

 Next we consider the following short exact commutative diagram\footnote{We 
pull back ${\cal X}\mapsto B$ by the mapping 
${\cal M}_{e_{II}}\mapsto B$ and abbreviate ${\cal M}_{e_{II}}\times_B{\cal X}$
by the same notation ${\cal X}$.},

\[
\hskip -.1in
\begin{array}{ccccc}
{\cal O}_{\cal X} & \mapsto & {\cal O}_{\cal X}(nD) & \mapsto& {\cal O}_{nD}(nD)\\
\Big\downarrow & &\Big\downarrow & & \Big\downarrow\\
{\cal O}_{\cal X}({\bf e}_{II})\otimes {\cal H}_{II} &\mapsto & 
{\cal O}_{\cal X}({\bf e}_{II}+nD) \otimes {\cal H}_{II} & \mapsto
& {\cal O}_{nD}({\bf e}_{II}+nD)\otimes {\cal H}_{II}\\
\Big\downarrow& & \Big\downarrow& & \Big\downarrow\\
{\cal O}_{{\bf e}_{II}}({\bf e}_{II})\otimes {\cal H}_{II} & 
\mapsto& {\cal O}_{{\bf e}_{II}}({\bf e}_{II}+nD) \otimes {\cal H}_{II}& \mapsto & {\cal O}_{{\bf e}_{II}\cap nD}({\bf e}_{II}+
nD)\otimes {\cal H}_{II}
\end{array}
\]

  This diagram is constructed from
 the twisted versions of the defining short exact sequences of the form
 $0\mapsto {\cal O}\mapsto {\cal O}({\bf D})\mapsto 
{\cal O}_{\bf D}({\bf D})\mapsto 0$ with
 ${\bf D}={\bf e}_{II}$, ${\bf e}_{II}$ and ${\bf e}_{II}|_{nD}$ 
for the columns and
 ${\bf D}=nD$, $nD$ and $nD|_{{\bf e}_{II}}$ for the rows, respectively. 
 By pushing these exact sequences forward along (the suitable restriction of) 
$\pi:{\cal M}_{e_{II}}\times_B
{\cal X}\mapsto {\cal M}_{e_{II}}$, we get the following commutative diagram
 of short exact sequences,

\[
\hskip -.3in
\begin{array}{ccccc}
{\cal R}^0\pi_{\ast}\bigl({\cal O}_{\cal X}\bigr) & \mapsto & 
{\cal R}^0\pi_{\ast}\bigl({\cal O}_{\cal X}(nD)\bigr)&\mapsto & 
{\cal R}^0\pi_{\ast}\bigl({\cal O}_{nD}(nD)\bigr)\\
\Big\downarrow & & \Big\downarrow& & \Big\downarrow\\
{\cal R}^0\pi_{\ast}\bigl({\cal O}_{\cal X}({\bf e}_{II})\bigr)\otimes {\cal H}_{II}
& \mapsto & 
{\cal R}^0\pi_{\ast}\bigl({\cal O}_{\cal X}({\bf e}_{II}+nD)\bigr)\otimes {\cal H}_{II}
 & \mapsto & 
{\cal R}^0\pi_{\ast}\bigl({\cal O}_{nD}({\bf e}_{II}+nD)\bigr)\otimes {\cal H}_{II}\\
\Big\downarrow & & \Big\downarrow& & \Big\downarrow\\
{\cal R}^0\pi_{\ast}\bigl({\cal O}_{{\bf e}_{II}}({\bf e}_{II})\bigr)\otimes {\cal H}_{II}
& \mapsto& {\cal R}^0\pi_{\ast}\bigl({\cal O}_{{\bf e}_{II}}({\bf e}_{II}+nD)\bigr)
\otimes {\cal H}_{II} & \mapsto
& {\cal R}^0\pi_{\ast}\bigl({\cal O}_{{\bf e}_{II}\cap nD}({\bf e}_{II}+nD)\bigr)
\otimes {\cal H}_{II}
\end{array}
\]

 As usual when we assume 
$e_{II}-c_1({\bf K}_{{\cal X}/B})$ is nef, the second derived image sheaf
 ${\cal R}^2\pi_{\ast}\bigl({\cal O}_{\cal X}({\bf e}_{II})\bigr)$ vanishes and
 we have the following sheaf surjection,

$$\hskip -.2in
{\cal R}^1\pi_{\ast}\bigl({\cal O}_{\cal X}({\bf e}_{II})
\bigr)\otimes {\cal H}_{II}
\mapsto {\cal R}^1\pi_{\ast}\bigl({\cal O}_{{\bf e}_{II}}({\bf e}_{II})
\bigr)\otimes {\cal H}_{II}
\mapsto {\cal R}^2\pi_{\ast}\bigl({\cal O}_{\cal X}\bigr)\mapsto 0.$$

 And this implies that ${\cal R}^1\pi_{\ast}\bigl({\cal O}_{{\bf e}_{II}}({\bf 
e}_{II}))\otimes {\cal H}_{II}$
 is mapped onto the locally free quotient sheaf 
${\cal R}^2\pi_{\ast}\bigl({\cal O}_{\cal X}\bigr)$ of rank $p_g$.

 On the other hand, we have the isomorphism 
 ${\cal R}^1\pi_{\ast}\bigl({\cal O}_{nD}(nD)\bigr)\cong 
{\cal R}^2\pi_{\ast}\bigl({\cal O}_{\cal X}\bigr)$, 
due to the vanishing of ${\cal R}^i\pi_{\ast}
\bigl({\cal O}_{\cal X}(nD)\bigr)$ for $i>0$ and $n\gg 0$. 

Thus we have the following commutative diagram of sheaves\footnote{These
two diagrams have overlapping blocks.},

\[
\hskip -.9in
\begin{array}{ccccccc}
{\cal R}^0\pi_{\ast}\bigl({\cal O}_{\cal X}({\bf e}_{II}+nD)\bigr)
\otimes {\cal H}_{II} 
& \mapsto & {\cal R}^0\pi_{\ast}\bigl({\cal O}_{nD}({\bf e}_{II}+nD)\bigr)
\otimes {\cal H}_{II} &\mapsto & 
{\cal R}^1\pi_{\ast}\bigl({\cal O}_{\cal X}({\bf e}_{II})\bigr)
\otimes {\cal H}_{II}& \mapsto & 0\\
\Big\downarrow & &\Big\downarrow & & \Big\downarrow & & \\
{\cal R}^0\pi_{\ast}\bigl({\cal O}_{{\bf e}_{II}}({\bf e}_{II}+nD)\bigr)
\otimes {\cal H}_{II} & \mapsto& 
{\cal R}^0\pi_{\ast}\bigl({\cal O}_{{\bf e}_{II}\cap nD}({\bf e}_{II}+nD)\bigr) 
\otimes {\cal H}_{II}& \mapsto& {\cal R}^1\pi_{\ast}\bigl({\cal O}_{{\bf e}_{II}}({\bf e}_{II})\bigr) \otimes {\cal H}_{II}& \mapsto 
& 0\\
\Big\downarrow & & \Big\downarrow& & \Big\downarrow & & \\
0 &\mapsto & {\cal R}^1\pi_{\ast}\bigl({\cal O}_{nD}(nD)\bigr)
& \stackrel{\cong}{\mapsto} & {\cal R}^2\pi_{\ast}\bigl({\cal O}_{\cal X}\bigr)
 & \mapsto & 0\\
& & \Big\downarrow & & \Big\downarrow & & \\
& & 0 & & 0 & &
\end{array}
\]

 Notice that ${\cal R}^0\pi_{\ast}\bigl({\cal O}_{nD}({\bf e}_{II}+nD)\bigr)$
 in the first row 
is nothing but the ${\bf W}_{e_{II}}$ in the datum of algebraic family
Kuranishi model of $e_{II}$. The second row is a part of a four-term exact
sequence regarding the fiberwise infinitesimal deformations and obstructions of
 ${\bf e}_{II}$.

 The above observation indicates that there exists
 a $4-$term exact sequence of obstruction vector bundles

$$\hskip -.7in
0\mapsto {\bf R}^0\pi_{\ast}\bigl({\cal O}_{nD}(nD)\bigr)\mapsto 
 {\bf R}^0\pi_{\ast}\bigl({\cal O}_{nD}({\bf e}_{II}+nD)\bigr)
\otimes {\bf H}_{II}
\mapsto {\bf R}^0\pi_{\ast}\bigl({\cal O}_{{\bf e}_{II}\cap nD}({\bf e}_{II}
+nD)\bigr)\otimes {\bf H}_{II}
\mapsto {\bf R}^2\pi_{\ast}\bigl({\cal O}_{\cal X}\bigr)\mapsto 0,$$

between the algebraic family obstruction bundle ${\bf W}_{e_{II}}$ 
of $e_{II}$ and the vector bundle 
${\bf R}^0\pi_{\ast}\bigl({\cal O}_{{\bf e}_{II}\cap nD}({\bf e}_{II}+nD)\bigr)
\otimes {\bf H}_{II}$
onto the infinitesimal obstructions 
${\bf R}^1\pi_{\ast}\bigl({\cal O}_{{\bf e}_{II}}({\bf e}_{II})\bigr)$. 
Their ranks differ by the geometric
 genus $p_g$ of the fiber algebraic surfaces of ${\cal X}\mapsto B$.

\medskip

As before let $ed$ denote the expected algebraic family Seiberg-Witten 
dimension of the
 type $II$ class $e_{II}$, $ed=
dim_{\bf C}B+p_g+{e_{II}^2-e_{II}\cdot c_1({\bf K}_{{\cal X}/B})
\over 2}$.

The following proposition shows that the virtual fundamental class
 of ${\cal M}_{e_{II}}$, $[{\cal M}_{e_{II}}]_{vir}$,
 appears naturally within the localized
 contribution of top Chern class of the bundle
 ${\bf R}^0\pi_{\ast}\bigl({\cal O}_{nD}(nD)
\bigr)\otimes
 {\bf H}_{II}\oplus {\bf W}'$ along ${\cal M}_{e_{II}}$.
It will play an essential role in our residual intersection theory 
approach in section \ref{subsection; stable}.

\begin{prop}\label{prop; local}
Let ${\bf W}_{e_{II}}'$ and ${\bf W}'$ be the vector bundles associated with 
the locally free sheaves ${\cal R}^0\pi_{\ast}\bigl({\cal O}_{{\bf e}_{II}\cap nD}
({\bf e}_{II}+nD)\bigr)\otimes {\cal H}_{II}$ and 
${\cal R}^0\pi_{\ast}\bigl({\cal O}_{{\bf e}_{II}\cap nD}
(\underline{C}+nD)\bigr)\otimes {\cal H}_{II}$ 
over ${\cal M}_{e_{II}}$, respectively. Then
 the localized contribution of top Chern class 
$$\hskip -.4in
\{c_{total}({\bf R}^0\pi_{\ast}\bigl({\cal O}_{nD}(nD)
\bigr)\otimes {\cal H}_{II}\oplus 
{\bf W}')\cap s({\cal M}_{e_{II}}, 
{\bf P}_B({\bf V}_{e_{II}}))\}_{ed-p_g}=
[{\cal M}_{e_{II}}]_{vir}\cap 
c_{p_g}({\cal R}^2\pi_{\ast}\bigl({\cal O}_{\cal X}\bigr))+\tilde{\eta}.$$
Over here $\tilde{\eta}$ is a cycle class which is a polynomial of 
the push-forward of monomials in $\underline{C}-e_{II}$, ${\bf e}_{II}$ and
 $nD$ along ${\bf e}_{II}\mapsto {\cal M}_{e_{II}}$.
\end{prop}

\noindent Proof of proposition \ref{prop; local}: We first recall that
 $[M_{e_{II}}]_{vir}=
\{c_{total}({\bf W}_{e_{II}}\otimes {\bf H}_{II}|_{{\cal M}_{e_{II}}})\cap 
s_{total}({\cal M}_{e_{II}}, 
{\bf P}({\bf V}_{e_{II}}))\}_{ed}\in {\cal A}_{\cdot}({\cal M}_{e_{II}})$, 
is the localized top Chern class of $\pi_{{\bf P}({\bf V}_{e_{II}})}^{\ast}
{\bf W}_{e_{II}}\otimes {\bf H}_{II}$. 
 On the other hand, we have observed from the above discussion
that $c_{total}({\bf W}_{e_{II}}')$ is the cap product of
 $c_{total}({\bf W}_{e_{II}}\otimes 
{\bf H}_{II})$ and $c_{total}({\bf R}^2\pi_{\ast}{\cal O}_{\cal X})$.

 So by capping with $s_{total}({\cal M}_{e_{II}}, {\bf P}({\bf V}_{e_{II}}))$ 
and by taking the degree $ed-p_g$ term, we find

$$
\{c_{total}(\bigl({\bf R}^0\pi_{\ast}{\cal O}_{nD}(nD)\otimes {\bf H}_{II}\oplus 
{\bf W}_{II}'\bigr)|_{{\cal M}_{e_{II}}})\cap s_{total}({\cal M}_{e_{II}}, 
{\bf P}({\bf V}_{e_{II}}))\}_{ed-p_g}$$
$$=\{c_{total}({\bf W}_{e_{II}}\otimes {\bf H}_{II}|_{{\cal M}_{e_{II}}})
\cap c_{total}({\bf R}^2\pi_{\ast}{\cal O}_{\cal X})
\cap s_{total}({\cal M}_{e_{II}}, 
{\bf P}({\bf V}_{e_{II}}))\}_{ed-p_g}$$
$$=[{\cal M}_{e_{II}}]_{vir}\cap c_{p_g}({\bf R}^2\pi_{\ast}{\cal O}_{\cal X}),$$

 by using that crucial property that the pairing 
$\{c_{total}({\bf W}_{e_{II}}\otimes {\bf H}_{II}
|_{{\cal M}_{e_{II}}})\cap s_{total}({\cal M}_{e_{II}}, 
{\bf P}({\bf V}_{e_{II}})\}_{ed-k}=0$, and 
$c_{p_g+k}({\bf R}^2\pi_{\ast}{\cal O}_{\cal X})=0$ for all $k>0$. 

Then the equality of our proposition follows from applying 
lemma \ref{lemm; factor}.
$\Box$

\medskip

\begin{rem}
 If the formal excess base dimension $febd(e_{II}, {\cal X}/B)=0$, then
 the expected dimension $ed=dim_{\bf C}B+{e_{II}^2-e_{II}\cdot 
c_1({\bf K}_{{\cal X}/B})\over 2}$. Then the identity in the above 
proposition should be replaced by 

$$\hskip -.4in
\{c_{total}({\bf R}^0\pi_{\ast}\bigl({\cal O}_{nD}(nD)
\bigr)\otimes {\cal H}_{II}\oplus 
{\bf W}')\cap s({\cal M}_{e_{II}}, 
{\bf P}_B({\bf V}_{e_{II}}))\}_{ed}=
[{\cal M}_{e_{II}}]_{vir}+\tilde{\eta}.$$
\end{rem}

\medskip

\section{The Blowup Construction of Algebraic Family Seiberg-Witten 
Invariants}\label{section; blowup}

\bigskip

 In this section we discuss the blowup construction of the
 family Seiberg-Witten invariant 
${\cal AFSW}_{{\cal X}\mapsto B}(1, \underline{C})$
with respect to a finite collection of
 type $II$ exceptional classes
 $e_{II; 1}$, $e_{II; 2}$, $\cdots$, $e_{II; p}$,
 generalizing the blowup and residual intersection theory 
construction of type $I$ exceptional
classes in [Liu6].  One major difference between the theories of
 type $I$ and type $II$ exceptional
classes is that for a type $II$ exceptional class $e_{II}$, the family
 moduli space of $e_{II}$, ${\cal M}_{e_{II}}=Z(s_{II})$ may
 not be regular and
 the cycle class $[{\cal M}_{e_{II}}]_{vir}={\bf Z}(s_{II})$ is
typically not equal to $[Z(s_{II})]$. 
 Actually for the canonical algebraic family Kuranishi model
 of a type $I$ exceptional class $e_i$ we \footnote{see section 6.2. of [Liu5].} 
  have ${\bf V}_{e_i}\cong {\bf C}$, the constant line bundle over $M_n$ and 
 ${\bf P}({\bf C})\cong M_n$. Moreover, the family moduli space of $e_i$,
 the existence locus of $e_i$ over $M_n$, can be identified with the
 closure of the admissible stratum $Y(\Gamma_{e_i})$ of the fan-like
 admissible graph $\Gamma_{e_i}$ (see the graph on page \pageref{fanlike}).

As $Y(\Gamma_{e_i})$ is smooth of 
 the expected dimension $dim_{\bf C}M_n+(e_i^2+1)$, 
$[Y(\Gamma_{e_i})]$ represents the fundamental
 class of the family moduli space ${\cal M}_{e_i}$.

Let ${\cal M}_{e_{II; i}}$ be the family moduli space of $e_{II; i}$ and 
 let $\pi_i:{\cal M}_{e_{II; i}}\mapsto B$ be the canonical projection into $B$.
 Over the locus $\pi_i({\cal M}_{e_{II; i}})\subset B$ the class $e_{II; i}$
 becomes effective and over their intersection 
$\cap_{1\leq i\leq p}\pi_i({\cal M}_{e_{II; i}})\subset B$ all the type $II$
 exceptional classes
 $e_{II; i}$ become effective simultaneously. 

\begin{defin}\label{defin; intersect}
Define $\cap_{1\leq i\leq p}\pi_i({\cal M}_{e_{II; i}})\subset B$ to be the locus
 of co-existence of $e_{II; 1}$, $e_{II; 2}$, $\cdots$, $e_{II; p}$.
Define ${\cal M}_{e_{II; 1}, \cdots, e_{II; p}}=
\times_B^{1\leq i\leq p}{\cal M}_{e_{II; i}}$ to be the moduli space
 of co-existence of the classes $e_{II; 1}$, $\cdots, e_{II; p}$.
\end{defin}

 Ideally we may expect ${\cal M}_{e_{II; i}}$ to be smooth of the
expected\footnote{Assuming $febd(e_{II}, {\cal X}/B)=p_g$.}
 dimension $dim_{\bf C}B+p_g+{e_{II}^2-c_1({\bf K}_{{\cal X}/B})\cdot
 e_{II}\over 2}$ and there exists a Zariski open and dense subset of 
${\cal M}_{e_{II; i}}$, called the ``interior'' of ${\cal M}_{e_{II; i}}$,
 parametrizing the irreducible curves representing $e_{II; i}$. 
Then ${\cal M}_{e_{II; i}}$
can be viewed as the natural compactification of its open and dense ``interior''.
 Under the idealistic assumption, we ``expect'' that there
exists a Zariski-dense open subset of 
$\times_B^{1\leq i\leq p}{\cal M}_{e_{II; i}}$ which parametrizes tuples of
 irreducible
universal type $II$ curves ${\bf e}_{II; i}$, $1\leq i\leq p$.

\bigskip

 In the real world the individual ${\cal M}_{e_{II; i}}$ may not be smooth,
 the intersection $\cap_{1\leq i\leq p}\pi_i({\cal M}_{e_{II; i}})$ or
the fiber product ${\cal M}_{e_{II; 1}, \cdots, e_{II; p}}=
\times_B^{1\leq i\leq p}{\cal M}_{e_{II; i}}$ is seldom regular.

\bigskip

 The basic philosophy of family Gromov-Taubes theory
 is to replace the objects ${\cal M}_{II; i}$ or
 ${\cal M}_{e_{II; 1}, \cdots, e_{II; p}}$ by the appropriated virtual 
 fundamental classes and interpret the enumeration of the invariants
 in term of intersection theory [F]. Thanks to the fact that all 
${\cal M}_{e_{II; i}}$ are compact, no complicated gluing construction
 is ever needed.

Because the numerical
 condition $e_{II; i}\cdot \underline{C}<0$ we impose on $\underline{C}$ and
 $e_{II; i}$, any effective representative
 of $\underline{C}$ over the ``interior'' of 
$\cap_{1\leq i\leq p}\pi_i({\cal M}_{e_{II; i}})$
 has to break off certain multiples of curves representing 
$e_{II; i}$, for each $1\leq i\leq p$. So we may write 
$\underline{C}=(\underline{C}-\sum e_{II; i})+\sum e_{II; i}$ formally.
 Thus we should be able to attach a family invariant of 
$\underline{C}-\sum_{1\leq i\leq p}e_{II; i}$ to (the
 virtual fundamental class) of $\times_B^{i\leq p}{\cal M}_{e_{II; i}}$
 using the geometric information
 of $e_{II; i}$, $1\leq i\leq p$ and express the localized contribution
 un-ambiguously.

\medskip

 In the following, we consider the following general question,

\medskip

\noindent {\bf Question}: Let $\underline{C}$ be an effective curve 
class over ${\cal X}\mapsto B$
 and let $e_{II; 1}, e_{II; 2}, e_{II; 3}, \cdots, e_{II; p}$ be $p$ distinct
 type $II$ exceptional classes over ${\cal X}\mapsto B$ such that
 $e_{II; i}\cdot \underline{C}<0$ for all $1\leq i\leq p$, while 
$e_{II; i}\cdot 
e_{II; j}\geq 0$ for $i\not=j$. What is the
 algebraic family Seiberg-Witten invariant attached to the
 moduli space of co-existence of $e_{II; i}$, $1\leq i\leq p$, 
$\times_B^{i\leq p}{\cal M}_{e_{II; i}}$? And what is the
 residual contribution of the algebraic family invariant of $\underline{C}$
 away from this moduli space of co-existence of $e_{II; i}$?

\medskip

 The resolution of the type $I$ 
analogue of the above question has been the 
 backbone of the proof of ``universality theorem'' [Liu6].

 Conceptually the residual contribution of the family invariant 
represents the contributions to the family
 invariants from curves in
 $\underline{C}$ within the family ${\cal X}\mapsto B$
 which are {\bf NOT} decomposed into a union of curves representing 
$\underline{C}-\sum_{1\leq i\leq p}e_{II; i}$ and $\sum_{1\leq i\leq p}e_{II; i}$,
 respectively.

\medskip

 There are a few important guidelines that we impose, based on
 the type $I$ theory, developed algebraically in [Liu6].

\medskip

\noindent {\bf Guideline 1}: We require that the localized (excess) contribution 
 of the family invariant of $\underline{C}$ 
along $\times_B^{1\leq i\leq p}{\cal M}_{e_{II; i}}$ to be proportional 
to the virtual fundamental class of ${\cal M}_{e_{II; 1}, \cdots, e_{II; p}}=
\times_B^{1\leq i\leq p}{\cal M}_{e_{II; i}}$ and to the virtual fundamental
class of ${\cal M}_{\underline{C}-\sum_{i\leq p}e_{II; i}}$.

In particular, when either $[{\cal M}_{e_{II; 1}, \cdots, e_{II; p}}]_{vir}$
 or $[{\cal M}_{\underline{C}-\sum_{i\leq p}e_{II; i}}]_{vir}$ vanishes,
 we require the desired localized contribution to vanish as well.

\medskip
\label{guide}

\noindent {\bf Guideline 2}: Because type $II$ exceptional classes can
 behave badly in comparison with their type $I$ siblings, 
the process of identifying
 the family invariant attached to $\underline{C}-\sum_{1\leq i\leq p}e_{II; i}$
 may be more delicate than the theory of type $I$ exceptional classes. But
 we expect that the resulting family invariant (see theorem \ref{theo; degenerate}
 for details)
 can be reduced to the modified algebraic family 
 Seiberg-Witten
\footnote{Consult definition 13 and 14 of [Liu6] for details.} 
invariant ${\cal AFSW}_{M_{n+1}\times T(M)\mapsto
 M_n\times T(M)}(c_{total}(\tau_{\Gamma}), C-{\bf M}(E)E-
\sum_{e_i\cdot (C-{\bf M}(E)E)<0}e_i)$ when $\underline{C}=C-{\bf M}(E)E$ and
 $e_{II; i}$ are reduced\footnote{The type $II$ curves satisfying condition 
 $febd(e_{II; i}, {\cal X}/B)=p_g$
 have different dimension formulae from the type $I$ curves'. This 
$p_g$ dimension shift introduces an additional $c_{p_g}({\cal R}^2\pi_{\ast}
{\cal O}_{\cal X})$ insertion for each type $II$ class.}
 to some collections of the type $I$ classes of the universal family 
$M_{n+1}\mapsto M_n$.

\medskip

\medskip

\noindent {\bf Guideline 3}: The localized (excess) contribution of the family 
invariant of $\underline{C}$ along $\times_B^{i\leq p} {\cal M}_{e_{II; i}}$ 
has to be independent to the algebraic family
Kuranishi models chosen for $\underline{C}$, $e_{II; i}$, $1\leq i\leq p$, etc.
 In particular, it is independent to $n\gg 0$ and the very ample divisor $D\subset
 {\cal X}$, etc, chosen to define the Kuranishi models.

\medskip

\label{guide4}
\noindent {\bf Guideline 4}: The construction we will provide should enable us to
 generalize to an inductive scheme involving more than one single collection of
 type $II$ exceptional classes.  Because exceptional curves can break up and
degenerate within a given family, instead of considering only 
a single collection of
 exceptional curves our scheme should work for a whole hierarchy of them.

\medskip

These few guidelines determine the localized (excess) contribution of the
 family invariant uniquely, as will be shown in theorem \ref{theo; degenerate}.

\medskip

 In the following subsection, some basic knowledge in intersection theory
[F] is recalled before we move on to the main theorem of the paper.

\medskip

\subsubsection{The Normal Cones and the Fiber Products}\label{subsubsection;
 cone}

\bigskip

  Let $X_i\mapsto B$, $1\leq i\leq n$, be $n$ purely $dim_{\bf C}X_i$ dimensional 
schemes over a smooth 
variety $B$ and let
 $Y_i\subset X_i$, $1\leq i\leq n$ be the closed sub-schemes of $X_i$
 defined by the zero loci of sections $s_i:X_i\mapsto E_i$ of vector 
bundles $E_i\mapsto X_i$.

 Consider the fiber products 
$\times_{B; i\leq n} Y_i\subset \times_{B; i\leq n}X_i$.
 In the sub-section we want to review the intersection product refined on
 $\times_{B; i\leq n} Y_i$, based on Fulton's general construction [F].

\medskip

  Recall that on page 132, definition 8.1.1. of [F], a refined product
 $x\cdot_f y$ is defined.  Let $f:X\mapsto Y$ with $Y$ non-singular and
let $p_X:X'\mapsto X$ and $P_Y:Y'\mapsto Y$ be morphisms of schemes.
Let $x\in {\cal A}_{\cdot}(X')$, $y\in {\cal A}_{\cdot}(Y')$. Then
 we can define $\gamma_f:X\mapsto X\times Y$ by $\gamma_{f}(t)=(t, f(t))$ 
and we have the following commutative diagram,

\[
\begin{array}{ccc}
X'\times_YY' & \mapsto & X'\times Y'\\
\Big\downarrow & & \Big\downarrow \\
X& \stackrel{\gamma_{f}}{\mapsto} & X\times Y
\end{array}
\]

\begin{defin} \label{defin; fiber}
Define $x\cdot_f y=\gamma_f^{!}(x\times y)$.
\end{defin}

Please consult page 132, proposition 8.1.1. for all the basic properties of
$\cdot_f$, including its associativity and commutativity, etc.

As usual, we take $[Y]_{vir}={\bf Z}(s_i)$.
 Our goal is to determine the virtual fundamental class of 
 $\times_B^{i\leq p}[Y_i]_{vir}$. 

\begin{defin}\label{defin; vir}
Define $[\times_B^{1\leq i\leq p}Y_i]_{vir}$ to be the Gysin pull-back
 $\Delta^{!}({\bf Z}(\oplus s_i))$ of $\oplus E_i\mapsto \times X_i$ by
 the diagonal morphism $\Delta:B\mapsto B^p$.
\end{defin}

 Based on mathematical induction, we
 may reduce to the $p=2$ case and prove the following proposition.

\medskip

\begin{prop}\label{prop; product}
Let $Y_1\subset X_1\mapsto B$ and $Y_2\subset X_2\mapsto B$ be closed 
sub-schemes over $B$ defined as the zero loci of the vector bundles 
 $E_i\mapsto X_i$. 
Then the virtual fundamental class of co-existence $[\times_B^{i\leq p} Y_i]_{vir}$
defined in definition \ref{defin; vir}
 is equal to $[Y_1]_{vir}\cdot_{id_B} [Y_2]_{vir}$.
\end{prop}

\noindent Proof of proposition \ref{prop; product}: Firstly consider the
 Cartesian product $Y_1\times Y_2\subset X_1\times X_2$. It is clear that
 $C_{Y_1\times Y_2}(X_1\times X_2)=C_{Y_1}X_1\times C_{Y_2}X_2$. The 
Cartesian products project naturally into $B\times B$ and  
the fiber product $Y_1\times_B Y_2$ or $X_1\times_B X_2$ can be viewed as the
 pull-back through $\Delta:B\mapsto B\times B$ of $Y_1\times Y_2\mapsto B\times B$
 or $X_1\times X_2\mapsto B\times B$.

  The virtual fundamental class of $Y_1\times Y_2$, $[Y_1\times Y_2]_{vir}$ 
is $$\{c_{total}(E_1\oplus E_2|_{Y_1\times Y_2})\cap s_{total}
(C_{Y_1\times Y_2}(X_1\times X_2))\}_{\sum dim_{\bf C}X_i-\sum rank_{\bf C}E_i}.$$

\medskip

 We know that $C_{Y_1\times Y_2}(X_1\times X_2)=C_{Y_1}X_1\times C_{Y_2}X_2$.

 From the following lemma, we can compute its total Segre class.

\begin{lemm}\label{lemm; pro}
Let $p_1:Y_1\times Y_2\mapsto Y_1$ and $p_2:Y_1\times Y_2\mapsto Y_2$ be the
 natural projections. Then the projections induce cones $p_1^{\ast}C_{Y_2}X_2$,
 $p_2^{\ast}C_{Y_1}X_1$ over $Y_1\times Y_2$, the normal cones of
 $Y_1\times Y_2\subset Y_1\times X_2$ and of $Y_1\times Y_2\subset
 X_1\times Y_2$.
We have the following identities on the total Segre classes,

$$s_{total}(C_{Y_1}X_1\times C_{Y_2}X_2)=s_{total}(p_2^{\ast}C_{Y_1}X_1)\cap
 s_{total}(p_1^{\ast}C_{Y_2}X_2)=s_{total}(C_{Y_1}X_1)\times
 s_{total}(C_{Y_2}X_2).$$
\end{lemm}

 The lemma is a generalization of the Whitney sum formula of
 vector bundles to normal cones. For completeness, we offer a simple proof here.

\medskip

\noindent Proof of lemma \ref{lemm; pro}: 

 If we have either $Y_1=X_1$ or $Y_2=X_2$, the above formula is a trivial identity.
 Let us assume $Y_i\not= X_i$, for $1\leq i\leq 2$. 
 We blow up $X_1, X_2$ along $Y_1, Y_2$, respectively and denote the resulting
 schemes by $\tilde{X}_1$, $\tilde{X}_2$, respectively. Let $D_1, D_2$ denote
the resulting exceptional divisors.

 Then we may blow up $\tilde{X}_1\times \tilde{X}_2$ along the 
codimension two $D_1\times D_2$ and
denote the resulting scheme $\tilde{X}_3$ with the exceptional divisor $D_3$.

  On the other hand, we may blow up $X_1\times X_2$ along $Y_1\times Y_2$ 
 directly and
denote the resulting scheme by $\widetilde{X_1\times X_2}$ and 
denote the exceptional divisor by
 $D$. Consider the dominated morphism $\tilde{X}_3\mapsto X_1\times X_2$,
 which maps $D_3$ onto $Y_1\times Y_2$. By the universal property of
scheme theoretical blowing up (proposition II.7.14 of [Ha]), 
the above map factors through
 $\tilde{X}_3\mapsto \widetilde{X_1\times X_2}$. We have the following commutative
 diagram, 

\[
\begin{array}{ccc}
 \tilde{X}_3 & \mapsto & \tilde{X}_1\times \tilde{X}_2\\
 \Big\downarrow & & \Big\downarrow \\
\widetilde{X_1\times X_2} & \mapsto & X_1\times X_2 \\
\end{array}
\]
 
 By using this commutative diagram, we may push-forward the total
Segre class of $C_{D_3}\tilde{X}_3$, $=\sum_{i\geq 0} c_1({\cal O}(-D_3))^i$, 
to $X_1\times X_2$ along two paths and the results must match. 
By using the birational invariance of
 the total Segre classes under proper birational push-forward (page 74, prop. 4.2
 of [F]), we conclude that

$$s_{total}(C_{Y_1\times Y_2}(X_1\times X_2))=s_{total}(p_2^{\ast}C_{Y_1}X_1)
\cap s_{total}(p_1^{\ast}C_{Y_2}X_2).$$
$\Box$

 By inserting the above 
identity into the defining equality of $[Y_1\times Y_2]_{vir}$,
 we find that $[Y_1\times Y_2]_{vir}=[Y_1]_{vir}\times [Y_2]_{vir}$.

 Moreover, we may take $X=Y=B$ and $f=id_B:B\mapsto B$, $X'=Y_1$ and
 $Y'=Y_2$, $x=[Y_1]_{vir}$ and $y=[Y_2]_{vir}$ in definition \ref{defin; fiber}.
 In this case $\gamma_f=\Delta:B\mapsto B\times B$ and the
 virtual fundamental class of the co-existence locus $[Y_1\times_B Y_2]_{vir}$
 is

$$\Delta^{!}([Y_1\times Y_2]_{vir})=
\Delta^{!}([Y_1]_{vir}\times [Y_2]_{vir})=[Y_1]_{vir}\times_{id_B} 
[Y_2]_{vir}.$$ $\Box$

 Because the refined intersection product $\cdot_f$ is associative, by
mathematical induction it is not hard to see that 
$[\times_B^{i\leq p} Y_i]_{vir}=\cdot_{id_B}^{1\leq i\leq p}[Y_i]_{vir}$. 

 We have the following simple lemma regarding their push-forwards into
 the global objects $\in{\cal A}_{\cdot}(B)$.

\begin{lemm}\label{lemm; cap}
Let $p_{Y_i}:Y_i\mapsto B$ denote the proper projection map 
from $Y_i$ to the smooth
 base space $B$.
The push-forward image of $[\times_B^{i\leq p}Y_i]_{vir}$ into 
${\cal A}_{\cdot}(B)$ (with $\cdot$ being the intersection product on $B$)
is the intersection product 
$p_{Y_1\ast}[Y_1]\cdot p_{Y_2\ast}
[Y_2]_{vir}\cdots ]_{vir}\cdot p_{Y_p\ast}[Y_p]_{vir}
\in {\cal A}_{\cdot}(B)$.
\end{lemm}

\noindent Proof of lemma \ref{lemm; cap}: Because $B$ is non-singular, 
the intersection product $\cdot$ makes 
${\cal A}^{\cdot}(B)={\cal A}_{dim_{\bf C}B-\cdot}(B)$ a 
commutative, graded ring with unit $[B]$.

 On the other hand, there is a commutative diagram,

\[
\begin{array}{ccc}
\times_B^{i\leq p}Y_i & \longrightarrow & \times_{i\leq p}Y_i\\
\Big\downarrow\vcenter{%
\rlap{$\scriptstyle{\mathrm{\times_{i\leq p}p_{Y_i}|_{\times_B^{i\leq p}Y_i}}}$}}&  & \Big\downarrow\vcenter{%
\rlap{$\scriptstyle{\mathrm{\times_{i\leq p}p_{Y_i}}}$}} \\
B& \stackrel{\Delta}{\longrightarrow} &  B^p
\end{array}
\]

\medskip

 By theorem 6.2. (a) on page 98 of [F], 
$$(\times_{i\leq p}p_{Y_i}|_{\times_B^{i\leq p}Y_i})_{\ast}\Delta^{!}=\Delta^{!}
(\times_{i\leq p}p_{Y_i})_{\ast}.$$

 Because $(\times_{i\leq p}p_{Y_i})_{\ast}(\times_{i\leq p}[Y_i]_{vir})=
\times_{i\leq p}p_{Y_i\ast}[Y_i]_{vir}$, the lemma follows from example 8.1.9. of
 [F]. $\Box$

\subsection{The Virtual Fundamental Class of 
${\cal M}_{e_{II; 1}, \cdots, e_{II; p}}$ and the Extension of
 Prop. \ref{prop; local}}\label{subsection; virtual}

\bigskip

We apply the discussion in subsection \ref{subsubsection; cone} to the
 concrete
situation of algebraic family Kuranishi models of $e_{II; i}$. 
Let $e_{II; i}, 1\leq i\leq p$ be $p$ distinct type $II$
 exceptional classes over ${\cal X}\mapsto B$. 
As in sub-section \ref{subsection; type2K}, we can take
 ${\cal V}_{II; i}$ and ${\cal W}_{II; i}$ to be 
${\cal R}^0\pi_{\ast}({\cal O}(nD)\otimes {\cal E}_{e_{II; i}})$ 
 and ${\cal R}^0\pi_{\ast}({\cal O}_{nD}(nD)\otimes {\cal E}_{e_{II; i}})$,
 respectively. Then $\Phi_{{\bf V}_{II; i}{\bf W}_{II; i}}:
{\bf V}_{II; i}\mapsto {\bf W}_{II; i}$ defines
 the algebraic family Kuranishi model for $e_{II; i}$.

 We take $X_i={\bf P}_B({\bf V}_{II; i})$ and $Y_i={\cal M}_{e_{II; i}}$.
 Then $[Y_i]_{vir}$ has to be defined to be the localized top Chern class of
$\pi_{{\bf P}({\bf V}_{II; i})}^{\ast}{\bf W}_{II; i}\otimes {\bf H}_{II; i}$,
 constructed in sub-section \ref{subsection; type2K}.

 By the general discussion in the preceding subsection 
\ref{subsubsection; cone}, prop. \ref{prop; product}, we may consider the
definition,

\begin{defin}\label{defin; coexist}
Define the virtual fundamental class 
$[{\cal M}_{e_{II; 1}, \cdots, e_{II; p}}]_{vir}$
of ${\cal M}_{e_{II; 1}, \cdots, e_{II; p}}=
\times_B^{1\leq i\leq p}{\cal M}_{e_{II; i}}$ to be
$\cdot_{id_B}^{1\leq i\leq p}[{\cal M}_{e_{II; i}}]_{vir}$.
\end{defin}

The following proposition is the natural 
extension of proposition \ref{prop; local} to
 the $p>1$ case.

\begin{prop}\label{prop; moreThanOne}
Let ${\cal H}_{II; i}$ be (the restriction of) the hyperplane invertible sheaf
 of ${\bf P}({\bf V}_{II; i})$ to ${\cal M}_{e_{II; 1}, \cdots, e_{II; p}}$.
Let ${\bf e}_{II; i}\mapsto {\cal M}_{e_{II; i}}$ be the universal curve
associated to each $e_{II; i}$ with $febd(e_{II; i}, {\cal X}/B)=p_g$.
 Let ${\bf W}'$ be the vector bundle associated to the locally free sheaf
 ${\cal R}^0\pi_{\ast}\bigl({\cal O}_{nD\cap \sum_{i\leq p}{\bf e}_{II; i}}(nD+
\underline{C})\bigr)\otimes_{i\leq p} {\cal H}_{II; i}$ over
 ${\cal M}_{e_{II; 1}, \cdots, e_{II; p}}$.
 Then the degree $dim_{\bf C}B+
\sum_{1\leq i\leq p}{e_{II; i}^2-c_1({\bf K}_{{\cal X}/B})\cdot
e_{II; i}\over 2}$ term of 

$$\hskip -.5in \{c_{total}(\oplus_{i\leq p} 
{\bf R}^0\pi_{\ast}\bigl({\cal O}_{nD}(nD)\bigr)\otimes
 {\cal H}_{II; i}\oplus 
{\bf W}'|_{{\cal M}_{e_{II; 1}, \cdots, e_{II; p}}})\cap 
s_{total}({\cal M}_{e_{II; 1}, \cdots, e_{II; p}}, \times_B^{1\leq i\leq p}
{\bf P}({\bf V}_{II; i}))\}$$ can be naturally expanded as 

$$[{\cal M}_{e_{II; 1}, \cdots, e_{II; p}}]_{vir}\cap 
c_{p_g}^p({\cal R}^2\pi_{\ast}{\cal O}_{\cal X})+\tilde{\eta}.$$

The class $\tilde{\eta}$ is a polynomial (in terms of $\cdot_{id_B}$) of
the push-forwards of algebraic expressions
 in terms of $\underline{C}-\sum_{i\in J}e_{II; i}$, $nD$ and the various 
$e_{II; i}$, where the index subset
$J$ runs through the subsets of $\{1, 2, \cdots, p\}$.
\end{prop}

 The proof of this proposition is based on mathematical 
induction and proposition
 \ref{prop; local}.

\noindent Sketch of the Proof: Firstly we notice that each
 ${\cal M}_{e_{II; i}}$ is of expected algebraic family dimension
 $dim_{\bf C}B+p_g+{e_{II; i}^2-c_1({\bf K}_{{\cal X}/B})\cdot e_{II; i}\over 2}$. 
 So ${\cal M}_{e_{II; 1}, \cdots, e_{II; p}}=\times_B^{1\leq i\leq p}
{\cal M}_{e_{II; i}}$ is of expected dimension 
$dim_{\bf C}B+p_g\cdot p+\sum_{i\leq p}
{e_{II; i}^2+c_1({\bf K}_{{\cal X}/B})\cdot e_{II; i}\over 2}$. So the
 degrees of the formula in the statement of the proposition match. When 
$p=1$, the above statement is reduced to proposition \ref{prop; local}.

 For general $p$, we consider the
Cartesian product $\times^{1\leq i\leq p}{\bf P}({\bf V}_{II; i})$
 and view the fiber product $\times_B^{1\leq i\leq p}{\bf P}({\bf V}_{II; i})$
 as its pull-back by the diagonal morphism $\Delta:B\mapsto B^p$. 

\medskip

  Consider the following sheaf short exact sequence,

$$\hskip -.5in
{\cal R}^0\pi_{\ast}\bigl({\cal O}_{\sum_{i\leq p-1}{\bf e}_{II; i}\cap nD}(
\underline{C}-e_{II; p}+nD)\bigr)\otimes^{i\leq p-1} {\cal H}_{II; i}
\mapsto {\cal R}^0\pi_{\ast}\bigl({\cal O}_{\sum_{i\leq p}{\bf e}_{II; i}\cap nD}(
\underline{C}+nD)\bigr)\otimes^{i\leq p} {\cal H}_{II; i}$$
$$\mapsto
{\cal R}^0\pi_{\ast}\bigl({\cal O}_{{\bf e}_{II; p}\cap nD}(\underline{C}+nD)\bigr)
\otimes^{i\leq p} {\cal H}_{II; i}.$$

 By our induction hypothesis for $p-1$, applying to the class $\underline{C}'=
\underline{C}-e_{II; p}$ and $p-1$ distinct exceptional classes 
$e_{II; 1}$, $e_{II; 2}$, $\cdots$, 
$e_{II; p-1}$, and the $p=1$ case (proposition \ref{prop; local}), applying to
 $\underline{C}'=\underline{C}$ and $e_{II; p}$, 
 using the computation of lemma \ref{lemm; pro}
 we can write the localized top Chern class as
$$\bigl([{\cal M}_{e_{II; 1}, \cdots, e_{II; p-1}}]_{vir}\cap 
c_{p_g}^{p-1}({\cal R}^2\pi_{\ast}{\cal O}_{\cal X})
+\tilde{\eta}_1
\bigr)\cdot_{id_B}\bigl([{\cal M}_{e_{II; p}}]_{vir}\cap 
c_{p_g}({\cal R}^2\pi_{\ast}{\cal O}_{\cal X})+\tilde{\eta}_2\bigr),$$

where $\tilde{\eta}_1$ and $\tilde{\eta}_2$ 
are $\cdot_{id_B}$-polynomial expressions of 
the push-forwards of $\underline{C}-e_{II; p}-\sum_{i\in I'}e_{II; i}$,
 $I'\subset \{1, \cdots, p-1\}$ and
 $\underline{C}-e_{II; p}$ and $nD$, $e_{II; i}$, etc. 
By a simple calculation, the conclusion
 follows. $\Box$ 

\medskip

\begin{rem}\label{rem; pg}
When $febd(e_{II; i}, {\cal X}/B)=0$, the term 
 $c_{p_g}^p({\cal R}^0\pi_{\ast}\bigl({\cal O}_{\cal X}\bigr))$ drops
 from the above formula in proposition \ref{prop; moreThanOne}.
\end{rem}

\subsection{The Stabilization of The Kuranishi Model of $\underline{C}$}
\label{subsection; stable}

\bigskip

 Let $(\Phi_{{\bf V}_{\underline{C}}{\bf W}_{\underline{C}}}, 
{\bf V}_{\underline{C}}, {\bf W}_{\underline{C}})$ be the algebraic family
 Kuranishi model of $\underline{C}$ defined by adopting some 
$nD$ with $n\gg 0$ and let 
$(\Phi_{{\bf V}_{II; i}{\bf W}_{II; i}}, {\bf V}_{II; i}, {\bf W}_{II; i})$
 be the algebraic family Kuranishi models of $e_{II; i}$, $1\leq i\leq p$,
 constructed following the recipe of subsection \ref{subsection; type2K}.

 Because the family moduli spaces ${\cal M}_{e_{II; i}}$ of the type $II$
 classes are embedded in
 ${\bf P}({\bf V}_{II; i})$ and not in $B$, we have to pull-back the
 algebraic family Kuranishi model of $\underline{C}$ from $B$ first.

\begin{lemm}\label{lemm; sub}
 Let ${\bf E}$ be a vector bundle over ${\cal T}_B({\cal X})$, then 
 ${\cal T}_B({\cal X})$
 can be identified with the subset of ${\bf P}({\bf E}\oplus {\bf C})$
 through the embedding induced by the bundle injection 
${\bf C}\mapsto {\bf E}\oplus {\bf C}$, and is the zero locus
 of the canonical section $s$ of 
$\pi_{{\bf P}({\bf E}\oplus {\bf C})}^{\ast}{\bf E}\otimes {\bf H}$ induced by 
${\bf E}\oplus {\bf C}\mapsto {\bf E}$.
\end{lemm}

\noindent Proof of lemma \ref{lemm; sub}: The cross-section $\sigma: 
{\cal T}_B({\cal X})\mapsto 
 {\bf P}({\bf E}\oplus {\bf C})$ induced by projectifying 
${\bf C}\mapsto {\bf E}\oplus {\bf C}$ is clearly isomorphic to 
${\cal T}_B({\cal X})$.
  On the other hand, the kernel of the bundle projection 
${\bf E}\oplus {\bf C}\mapsto {\bf E}$ is exactly the trivial sub-bundle ${\bf C}$.
 This implies that ${\bf H}^{\ast}\mapsto 
\pi_{{\bf P}({\bf E}\oplus {\bf C})}^{\ast}{\bf E}$ induced by
 ${\bf E}\oplus {\bf C}\mapsto {\bf E}$ is injective exactly
 off $\sigma({\cal T}_B({\cal X}))\subset {\bf P}({\bf E}\oplus {\bf C})$. So the
 canonical section $s$ of 
$\pi_{{\bf P}({\bf E}\oplus {\bf C})}^{\ast}{\bf E}\otimes {\bf H}$ vanishes
 exactly on $\sigma({\cal T}_B{\cal X})$ 
and it is easy to see that $s$ is a regular section
 of $\pi_{{\bf P}({\bf E}\oplus {\bf C})}^{\ast}{\bf E}\otimes {\bf H}$.
 $\Box$

\begin{rem}\label{rem; cutdown}
The cycle class of the zero locus $\sigma({\cal T}_B({\cal X}))$,
 $[\sigma({\cal T}_B({\cal X}))]\in 
{\cal A}_{\cdot}({\bf P}({\bf E}\oplus {\bf C}))$ is
 equal to 
$c_{top}(\pi_{{\bf P}({\bf E}\oplus {\bf C})}^{\ast}{\bf E}\otimes {\bf H})\cap
 [{\bf P}({\bf E}\oplus {\bf C})]$.
\end{rem}

\medskip

\begin{defin}\label{defin; Bp}
Denote the fundamental cycle class of the
zero cross-section $B\mapsto \times_B^p{\cal T}_B({\cal X})$ to be 
 $[B]_p$.
\end{defin} 

 The normal bundle of $[B]_p$ in $\times_B^p{\cal T}_B({\cal X})$ is isomorphic
 to ${\bf R}^1\pi_{\ast}{\cal O}_{\cal X}^{\oplus p}$.

  We replace $B$ by the auxiliary space 
 $B'=\times_B^{1\leq i\leq p}{\bf P}({\bf V}_{II; i}\oplus {\bf C})$
 and view the original $\times_B^p{\cal T}_B({\cal X}) \subset 
B'=\times_B^{1\leq i\leq p}
{\bf P}({\bf V}_{II; i}\oplus {\bf C})$
 as the regular 
zero locus of the canonical section of the auxiliary obstruction bundle
 $\oplus_{1\leq i\leq p}
\pi_{{\bf P}({\bf V}_{II; i}\oplus {\bf C})}^{\ast}{\bf V}_{II; i}\otimes
 {\bf H}_{II; i}$.

 Then by remark \ref{rem; cutdown} and definition \ref{defin; Bp}
we have to compensate by inserting both the
 top Chern class 
$c_{top}(\oplus_{1\leq i\leq p}\pi_{{\bf P}({\bf V}_{II; i}\oplus {\bf C})}^{\ast}
{\bf V}_{II; i}\otimes {\bf H}_{II; i})$ and the $[B]_p$
into the intersection pairing of the 
family invariant of $\underline{C}$.

\medskip

 Correspondingly for these $p$ distinct 
type $II$ classes $e_{II; i}$, we stabilize their algebraic family
 Kuranishi models by the trivial line bundle ${\bf C}$ and get
 $(\Phi_{{\bf V}_{II; i}{\bf W}_{II; i}}\oplus id_{\bf C}, 
 {\bf V}_{II; i}\oplus {\bf C}, {\bf W}_{II; i}\oplus {\bf C})$.

 By lemma \ref{lemm; stable}, these models are invariant under stabilizations.
 Then the push-forward image of the virtual 
fundamental class of ${\cal M}_{e_{II; i}}$
 into ${\cal A}_{\cdot}({\bf P}({\bf V}_{II; i}\oplus {\bf C}))$ is equal to

$$
 c_{top}(\pi_{{\bf P}({\bf V}_{II; i}\oplus {\bf C})}^{\ast}({\bf W}_{II; i}
\oplus {\bf C})\otimes {\bf H}_{II; i})\cap 
[{\bf P}({\bf V}_{II; i}\oplus {\bf C})]$$
$$=
c_{top}(\pi_{{\bf P}({\bf V}_{II; i}\oplus {\bf C})}^{\ast}{\bf W}_{II; i}\otimes
 {\bf H}_{II; i})\cap c_1({\bf H}_{II; i})\cap[{\bf P}({\bf V}_{II; i}\oplus 
{\bf C})]$$

$$=c_{top}(\pi_{{\bf P}({\bf V}_{II; i}\oplus {\bf C})}^{\ast}{\bf W}_{II; i}
\otimes
 {\bf H}_{II; i})\cap [{\bf P}({\bf V}_{II; i})],$$

 \footnote{${\bf P}({\bf V}_{II; i})$ can be viewed as the compactifying
 divisor at the infinity of ${\bf P}({\bf V}_{II; i}\oplus {\bf C})$.}
because 
$[{\bf P}({\bf V}_{II; i})]=c_1({\bf H}_{II; i})\cap[{\bf P}({\bf V}_{II; i}
\oplus {\bf C})]$.

\bigskip

 We introduce the following short-hand notation which will be used 
frequently later.

\begin{defin}\label{defin; obstruction}
Define ${\bf U}_{e_{II; 1}, \cdots, e_{II; p}}$ to be

$$\oplus_{1\leq i\leq p}\pi_{\times_B^{1\leq i\leq p}
 {\bf P}({\bf V}_{II; i}\oplus
 {\bf C})}^{\ast}({\bf V}_{II; i})\otimes {\bf H}_{II; i}.$$
\end{defin}

   The family moduli space ${\cal M}_{\underline{C}}$ is 
embedded in ${\bf P}({\bf V}_{\underline{C}})$ as a projectified abelian
 cone and it is the zero locus $Z(s_{\underline{C}})$ 
of a section $s_{\underline{C}}$ of 
$\pi_{{\bf P}({\bf V}_{\underline{C}})}^{\ast}{\bf W}_{\underline{C}}\otimes 
 {\bf H}$. The algebraic family Seiberg-Witten invariant of $\underline{C}$
is defined to be the integral of the top intersection pairing of
 $c_{top}(\pi_{{\bf P}({\bf V}_{\underline{C}})}^{\ast}
{\bf W}_{\underline{C}}\otimes 
 {\bf H})$ capping with a suitable power of $c_1({\bf H})$.

 We pull back the algebraic family Kuranishi model 
$(\Phi_{{\bf V}_{\underline{C}}{\bf W}_{\underline{C}}}, {\bf V}_{\underline{C}}, 
{\bf W}_{\underline{C}})$ from $B$ to $\times_B^{1\leq i\leq p}
{\bf P}({\bf V}_{II; i}\oplus {\bf C})$. To simplify our notation, we 
 skip the pull-back notation and denote the datum of the
 Kuranishi models by the same symbols.

 By remark \ref{rem; cutdown}, we have to extend the obstruction bundle from
$\pi_{{\bf P}({\bf V}_{\underline{C}})}^{\ast}
{\bf W}_{\underline{C}}\otimes 
 {\bf H}$ to $\pi_{{\bf P}({\bf V}_{\underline{C}})}^{\ast}
{\bf W}_{\underline{C}}\otimes 
 {\bf H}\oplus {\bf U}_{e_{II; 1}, \cdots, e_{II; p}}$, or 
equivalently to
$\pi_{{\bf P}({\bf V}_{\underline{C}})}^{\ast}
{\bf W}_{\underline{C}}\otimes 
 {\bf H}\otimes_{i\leq p}{\bf H}_{II; i}
\oplus {\bf U}_{e_{II; 1}, e_{II; 2}, \cdots, e_{II; p}}$. And then insert
 $[B]_p$ into the intersection pairing.

\begin{rem}\label{rem; twist}
 The above twisting of $\pi_{{\bf P}({\bf V}_{\underline{C}})}^{\ast}
{\bf W}_{\underline{C}}\otimes {\bf H}$ 
by $\otimes_{i\leq p}{\bf H}_{II; i}$ does not
affect the family invariant because the embedding $\sigma:{\cal T}_B({\cal X}) 
\mapsto {\bf P}({\bf V}_{II; i}\oplus {\bf C})$ defined by remark
 \ref{rem; cutdown} is totally
disjoint from the smooth divisor at infinity $\cong {\bf P}({\bf V}_{II; i})$ and
the line bundle ${\bf H}_{II; i}$ is trivialized over 
 $\sigma({\cal T}_B({\cal X}))$.
\end{rem}

 On the other hand, the moduli space of 
co-existence of $e_{II; 1}, \cdots, e_{II; p}$,
 ${\cal M}_{e_{II; 1}, \cdots, e_{II; p}}$, is
 a closed sub-scheme of the auxiliary base space
$B'=\times_B^{1\leq i\leq p}{\bf P}({\bf V}_{II; i}\oplus {\bf C})$. 
 So ${\cal M}_{\underline{C}}
\times_{B'}
 {\cal M}_{e_{II; 1}, e_{II; 2}, \cdots, e_{II; p}}$ is also a closed sub-scheme
 of the projective bundle $X={\bf P}({\bf V}_{\underline{C}})$.  By 
section II.7, page 160-161 of [Ha], we may blow up
 $X={\bf P}_{B'}({\bf V}_{\underline{C}})$ along
 $Z(s_{\underline{C}})
\times_{B'}
 {\cal M}_{e_{II; 1}, \cdots, e_{II; p}}={\cal M}_{\underline{C}}
\times_{B'}{\cal M}_{e_{II; 1}, \cdots, e_{II; p}}$ and make it into a divisor 
 ${\bf D}$ of the
 blown up scheme $\tilde{X}$.  And the direct 
application of the residual intersection
 formula implies that we may rewrite

$$
c_{top}(\pi_{\tilde{X}}^{\ast}{\bf W}_{\underline{C}}\otimes {\bf H}
\otimes_{j\leq p}
 {\bf H}_{II; j})=
c_{top}(\pi_{\tilde{X}}^{\ast}{\bf W}_{\underline{C}}\otimes {\bf H}
\otimes_{j\leq p} {\bf H}_{II; j}\otimes 
{\cal O}(-{\bf D}))$$
$$+\sum_{1\leq i\leq rank_{\bf C}{\bf W}_{\underline{C}}}(-1)^{i-1} c_{rank_{\bf C}
{\bf W}_{\underline{C}}-i}(
\pi_{\tilde{X}}^{\ast}{\bf W}_{\underline{C}}\otimes 
{\bf H}\otimes_{j\leq p}{\bf H}_{II; j})\cap
{\bf D}^{i-1}[{\bf D}].$$

 And the push-forward of 
$$\hskip -.2in
\sum_{1\leq i\leq rank_{\bf C}{\bf W}_{\underline{C}}}(-1)^{i-1} c_{rank_{\bf C}
{\bf W}_{\underline{C}}-i}(\pi_{\tilde{X}}^{\ast}{\bf W}_{\underline{C}}\otimes 
{\bf H}\otimes_{j\leq p}{\bf H}_{II; j}\oplus {\bf U}_{e_{II; 1}, \cdots,
 e_{II; p}})$$
$$\cap
{\bf D}^{i-1}[{\bf D}]\cap [B]_p\cap 
c_1({\bf H})^{dim_{\bf C}B+p_g-q+{\underline{C}^2-
c_1({\bf K}_{{\cal X}/B})\cdot \underline{C}\over 2}}$$
 
 is the localized contribution of the algebraic family Seiberg-Witten invariant
 along ${\cal M}_{e_{II; 1}, \cdots, e_{II; p}}$.

\medskip

 It is vital to understand:

\bigskip

\noindent {\bf Question}: Is the localized contribution of top Chern class
 the correct ``invariant'' associated to the collection of type $II$ 
classes $e_{II; 1}$, $\cdots$, $e_{II; p}$? In other words, is the
 localized contribution of top Chern classes constructed above 
invariant under deformations of 
the family ${\cal X}\mapsto B$ or of the Kuranishi models?

\bigskip

 When the exceptional classes are of type $I$, it has been shown in
 [Liu5], [Liu6] that the localized contribution 
of top Chern class can be
identified with certain mixed family invariant of $C-{\bf M}(E)E-\sum_{e_i\cdot
(C-{\bf M}(E)E)<0}e_i$ and consequentially are known to be topological. 
On the other hand, the theory for type $II$ curves
 is more delicate than their type $I$ counterpart as the 
naively chosen localized
contribution may not always be an invariant.

\bigskip

 To understand what may go wrong,
  one may consider some deformation of the 
datum $\Phi_{{\bf V}_{II}{\bf W}_{II}}:
{\bf V}_{II}\mapsto {\bf W}_{II}$. In order for it
 the correct choice, the localized contribution of
 top Chern class has to be 
independent to the deformations.  Consider the idealistic
 situation that
under the one parameter family of 
deformations, the family moduli space ${\cal M}_{e_{II}}$ of a single 
$e_{II}$ (i.e. $p=1$) 
is deformed into the whole space ${\bf P}({\bf V}_{e_{II}}\oplus {\bf C})$. 
After such
a degeneration\footnote{Geometrically this corresponds to thickening the
 family moduli space of $e_{II}$ to the whole space, which can be achieved by
 multiplying the map $\Phi_{{\bf V}_{e_{II}}{\bf W}_{e_{II}}}$ by 
a constant $t$ and shrink $t$ to zero.}, the family moduli space of $e_{II}$ 
is apparently not of the expected dimension. It is easy to see that the localized
contribution of top Chern class along such a degenerated family moduli space is
nothing but the whole localized top Chern class of ${\bf H}\otimes 
\pi_X^{\ast}{\bf W}_{\underline{C}}$ along ${\cal M}_{\underline{C}}$.

On the other hand, for the original well-behaved ${\cal M}_{e_{II}}$ the localized
 contribution of top Chern class is usually not equal to the whole
 localized top Chern class of ${\bf H}\otimes 
\pi_X^{\ast}{\bf W}_{\underline{C}}$ along ${\cal M}_{\underline{C}}$. 
Therefore we observe in this hypothetical
example that the localized contribution of top Chern classes may {\bf not}
 be invariant to the degenerations.

\medskip

 We also realize from this example that the non-invariant nature of the
localized contribution of top Chern classes is due to the non-invariant 
nature of the family moduli space ${\cal M}_{e_{II}}$! 

\bigskip

\label{special}
\noindent {\bf Special Assumption}: In the following, we assume that

\noindent 
${\cal M}_{\underline{C}-\sum e_{II; i}}\times_B{\cal M}_{e_{II; 1}, \cdots
, e_{II; p}}\hookrightarrow {\cal M}_{\underline{C}}
\times_B{\cal M}_{e_{II; 1}, \cdots, e_{II; p}}$ induced by
 adjoining the union of type $II$ curves in $\sum_{i\leq p} e_{II; i}$ to a curve
 in $\underline{C}-\sum_{i\leq p} e_{II; i}$ has been an isomorphism.
 
\medskip

 This assumption is the analogue of the simplifying assumption of
theorem 4. in [Liu5].

 The following is the main theorem in this paper,

\begin{theo}\label{theo; degenerate}
 Given the localized contribution of top Chern class
$$\hskip -1in
\sum_{1\leq i\leq rank_{\bf C}{\bf W}_{\underline{C}}}(-1)^{i-1} c_{rank_{\bf C}
{\bf W}_{\underline{C}}-i}(\pi_{\tilde{X}}^{\ast}{\bf W}_{\underline{C}}\otimes 
{\bf H}\otimes_{j\leq p}{\bf H}_{II; j}\oplus {\bf U}_{e_{II; 1}, \cdots,
 e_{II; p}})$$
$$\cap
{\bf D}^{i-1}[{\bf D}]\cap [B]_p\cap 
c_1({\bf H})^{dim_{\bf C}B+p_g-q+{\underline{C}^2-
c_1({\bf K}_{{\cal X}/B})\cdot \underline{C}\over 2}},$$

then under the above {\bf special assumption},
it can be expanded into an algebraic expression of cycle classes.
 Among the various terms of the expansion, there is a dominating term
proportional to the virtual fundamental class 
$[{\cal M}_{e_{II; 1}, \cdots, e_{II; p}}]_{vir}$, and is identified to be
 
$${\cal AFSW}_{{\cal X}\mapsto B}((\times_B^{i\leq p}
 \pi_i)_{\ast}\bigl([{\cal M}_{e_{II; 1}, \cdots, e_{II; p}}]_{vir}
\cap [B]_p\bigr)\cap c_{p_g}^p({\cal R}^2\pi_{\ast}{\cal O}_{\cal X})
 \cap \tau, \underline{C}-\sum_{1\leq i\leq p}e_{II; i}).$$

It satisfies the
 crucial invariant properties,

\medskip

(i). The class $\tau\in 
{\cal A}_{\cdot}({\cal M}_{e_{II; 1}, \cdots, e_{II; p}})$
 is independent to $nD$ and is depending on $\underline{C}$ and
 $e_{II; 1}$, $\cdots$, $e_{II; p}$ only. 

\medskip

(ii). When $febd(e_{II; i}, {\cal X}/B)=p_g$, $1\leq i\leq p$, and 
${\cal M}_{e_{II; 1}, \cdots, e_{II; p}}$ 
is smooth of its expected dimension $dim_{\bf C}B+p\cdot p_g+\sum {e_{II; i}^2-
c_1({\bf K}_{\cal X}/B)\cdot e_{II; i}\over 2}$,
 the smooth cycle $[{\cal M}_{e_{II; 1}, \cdots, e_{II; p}}]$ 
coincide with the virtual fundamental
class $[{\cal M}_{e_{II; 1}, \cdots, e_{II; p}}]_{vir}$ 
and the above mixed invariant can be re-expressed as
$${\cal AFSW}_{{\cal X}\times_B{\cal M}_{e_{II; 1}, \cdots, e_{II; p}}\mapsto 
{\cal M}_{e_{II; 1}, \cdots, e_{II; p}}}([B]_p\cap \tau, 
\underline{C}-\sum_{1\leq i\leq p}e_{II; i}).$$

\medskip

(iii). The above expression satisfies guidelines 1-3 listed beginning from
 page \pageref{guide}.
\end{theo}

\noindent Proof of theorem \ref{theo; degenerate}: After we push-forward
along ${\bf D}\mapsto {\cal M}_{e_{II; 1}, \cdots, e_{II; p}}$,
the localized contribution of top Chern classes of 
$\pi_X^{\ast}{\bf W}_{\underline{C}}\otimes {\bf H}\otimes_{i\leq p}
 {\bf H}_{II; i}$ 
along ${\cal M}_{e_{II; 1}, \cdots, e_{II; p}}$ is equal to 
 $\{c_{total}(\pi_X^{\ast}{\bf W}_{\underline{C}}\otimes {\bf H}
\otimes_{i\leq p} {\bf H}_{II; i}\oplus 
{\bf U}_{e_{II; 1}, e_{II; 2}, \cdots, e_{II; p}})
\cap s_{total}({\cal M}_{\underline{C}}\times_B {\cal M}_{e_{II; 1}, \cdots
, e_{II; p}}, {\bf P}({\bf V}_{\underline{C}})\times_B
B'))\cap [B]_p\}_{dim_{\bf C}B+p_g-q+{\underline{C}^2-\underline{C}\cdot 
c_1({\bf K}_{{\cal X}/B})\over 2}}$.

Assuming that ${\cal M}_{\underline{C}-\sum_{i\leq p}e_{II; i}}\times_B 
{\cal M}_{e_{II; 1}, \cdots, e_{II; p}}={\cal M}_{\underline{C}}\times_B 
{\cal M}_{e_{II; 1}, \cdots, e_{II; p}}$, the Segre class of the normal cone of 
${\cal M}_{\underline{C}}\times_B{\cal M}_{e_{II; 1}, \cdots, e_{II; p}}$
 is the same as the Segre class of 
${\cal M}_{\underline{C}-\sum_{i\leq p}e_{II; i}}
\times_B{\cal M}_{e_{II; 1}, \cdots, e_{II; p}}\subset 
{\bf P}({\bf V}_{\underline{C}})$.

\noindent Step I: Recall the following short exact sequence\footnote{Check
 the statement of 
proposition \ref{prop; moreThanOne} for the definition of ${\bf W}'$.}
in subsection \ref{subsection; virtual}, 

$$0\mapsto {\bf W}_{\underline{C}-\sum_{i\leq p}e_{II; i}}\mapsto 
{\bf W}_{\underline{C}}\otimes_{i\leq p}{\bf H}_{II; i}\mapsto {\bf W}'\mapsto 0.$$

 By applying lemma \ref{lemm; pro} and the discussion in subsection 
\ref{subsubsection; cone}, we can view the above localized contribution of
 top Chern class as the Gysin pull-back $\Delta^{!}$ (with 
$\Delta:B\mapsto B\times B$) from ${\cal A}_{\cdot}(
{\cal M}_{\underline{C}-\sum e_{II; i}}
\times {\cal M}_{e_{II; 1}, \cdots, e_{II; p}})$.
 We may rewrite the original
 localized contribution of top Chern class as the $\Delta^{!}$ pull-back of

 $$\{c_{total}(\pi_X^{\ast}{\bf W}_{\underline{C}-\sum e_{II; i}}\otimes {\bf H})
\cap s_{total}({\cal M}_{\underline{C}-\sum_{i\leq p}e_{II; i}}, 
{\bf P}({\bf V}_{\underline{C}-\sum e_{II; i}}))$$

$$\times c_{total}( 
{\bf W}'\otimes {\bf H})\cap s_{total}({\cal M}_{e_{II; 1}, 
e_{II; 2}, \cdots, e_{II; p}}, B')\cap c_{total}({\bf U}_{e_{II; 1}, \cdots,
 e_{II; p}})$$
$$\cap s_{total}({\bf V}'\otimes {\bf H})\cap [B]_p
\}_{2dim_{\bf C}B+p_g-q+{\underline{C}^2-
c_1({\bf K}_{{\cal X}/B})\cdot \underline{C}\over 2}}.$$

 We have used the short exact \footnote{Check the commutative 
diagram in proposition \ref{prop; compare} for the $p=1$ version.}
 sequence $0\mapsto 
{\bf V}_{\underline{C}-\sum e_{II; i}}\mapsto
 {\bf V}\otimes_{i\leq p} {\bf H}_{II; i}\mapsto {\bf V}'\mapsto 0$ over
 ${\cal M}_{e_{II; 1}, \cdots, e_{II; p}}$ and the following identity

$$
s_{total}({\cal M}_{\underline{C}-\sum e_{II; i}}\times_{B'}
{\cal M}_{e_{II; 1}, \cdots, e_{II; p}}, {\bf P}({\bf V}_{\underline{C}}))$$
$$=s_{total}({\cal M}_{\underline{C}-\sum e_{II; i}}\times_{B'}
{\cal M}_{e_{II; 1}, \cdots, e_{II; p}}, {\bf P}({\bf V}_{\underline{C}
-\sum e_{II; i}}))\cap s_{total}({\bf V}'\otimes {\bf H}).$$

Based on the following identity

$$\hskip -.3in 
dim_{\bf C}B-q+p_g+{(\underline{C}-\sum_{i\leq p}
 e_{II; i})^2-c_1({\bf K}_{{\cal X}/B})
\cdot (\underline{C}-\sum_{i\leq p} e_{II; i})\over 2}+\sum_{i\leq p} 
{ e_{II; i}^2-c_1({\bf K}_{{\cal X}/B})
\cdot e_{II; i}\over 2}$$
$$=dim_{\bf C}B-q+p_g+{\underline{C}^2-
c_1({\bf K}_{{\cal X}/B})
\cdot \underline{C}\over 2}+\{\sum_{i\leq p}
(-\underline{C}\cdot e_{II; i}+e_{II; i}^2)+
\sum_{1\leq i<j\leq p}e_{II; i}\cdot e_{II; j} \},$$

 on the family dimensions, we may set:
 
\noindent $a_1=dim_{\bf C}B-q+p_g+{(\underline{C}-\sum e_{II; i})^2-
c_1({\bf K}_{{\cal X}/B})
\cdot (\underline{C}-\sum e_{II; i})\over 2}$, $a_2=dim_{\bf C}B+\sum 
{e_{II; i}^2-c_1({\bf K}_{{\cal X}/B})
\cdot e_{II; i}\over 2}$ and $a_3=\sum(-\underline{C}\cdot e_{II; i}+e_{II; i}^2)+
\sum_{1\leq i<j\leq p}e_{II; i}\cdot e_{II; j}$ and focus on the term
\footnote{Assume that $a_3$ satisfies $a_3-p\cdot q(M)\geq 0$.
 Also notice that we have the identity $a_1+a_2-a_3=2dim_{\bf C}B-q(M)+p_g
+{\underline{C}^2-c_1({\bf K}_{{\cal X}/B})\cdot \underline{C}\over 2}$.}

$$g_{a_1, a_2, a_3}=
\{c_{total}(\pi_X^{\ast}{\bf W}_{\underline{C}-\sum e_{II; i}}\otimes {\bf H})
\cap s_{total}({\cal M}_{\underline{C}-\sum_{i\leq p}e_{II; i}}, 
{\bf P}({\bf V}_{\underline{C}-\sum e_{II; i}}))\}_{a_1}$$

$$\times \{c_{total}(\oplus_{i\leq p}{\bf R}^0\pi_{\ast}\bigl({
\cal O}_{nD}(nD)\bigr)\otimes {\bf H}_{II; i}\oplus 
{\bf W}'\otimes {\bf H})\cap s_{total}({\cal M}_{e_{II; 1}, 
e_{II; 2}, \cdots, e_{II; p}}, B')\cap_{i\leq p}c_1({\bf H}_{II; i})
\}_{a_2}$$
$$\cap \{c_{a_3-p\cdot q(M)}({\bf U}_{e_{II; 1}, \cdots,
 e_{II; p}}-\oplus_{i\leq p}({\bf C}\oplus {\bf R}^0\pi_{\ast}\bigl({
\cal O}_{nD}(nD)\bigr))\otimes {\bf H}_{II; i}-{\bf V}'\otimes {\bf H})
\cap [B]_p\}.$$

 The reason to pick this particular combination
 will be clear momentarily.

\bigskip

\noindent (i). By using $\cap_{i\leq p} c_1({\bf H}_{II; i})\cap
 [B']=[\times_B^{i\leq p}{\bf P}({\bf V}_{II; i})]$, 
the expression 
$\{c_{total}(\pi_X^{\ast}{\bf W}_{\underline{C}-\sum e_{II; i}}\otimes {\bf H})
\cap s_{total}({\cal M}_{\underline{C}-\sum_{i\leq p}e_{II; i}}, 
{\bf P}({\bf V}_{\underline{C}-\sum e_{II; i}}))\}_{a_1}$ is nothing but the
 virtual fundamental class of ${\cal M}_{\underline{C}-\sum e_{II; i}}$,
expressed as the localized top Chern class of 
 $\pi_X^{\ast}{\bf W}_{\underline{C}-\sum e_{II; i}}\otimes {\bf H}$.

\medskip

\noindent (ii). By applying proposition \ref{prop; moreThanOne} to
\footnote{Notice that we have added $\oplus_{i\leq p}
 {\bf C}\otimes {\bf H}_{II; i}$ to our bundle here. It is to match
 with the additional ${\bf H}_{II; i}$ factor in the
 stabilized obstruction bundle of $e_{II; i}$,
 $(\pi_{{\bf P}({\bf V}_{II; i}\oplus {\bf C})}^{\ast}
{\bf W}_{II; i}\oplus {\bf C})\otimes {\bf H}_{II; i}$.}
$\{c_{total}((\oplus_{i\leq p}{\bf R}^0\pi_{\ast}\bigl({
\cal O}_{nD}(nD)\bigr)\oplus {\bf C})\otimes {\bf H}_{II; i}\oplus 
{\bf W}'\otimes {\bf H})\cap s_{total}({\cal M}_{e_{II; 1}, 
e_{II; 2}, \cdots, e_{II; p}}, B')\}_{a_2}$, it can be expanded into 
an algebraic expression of cycle classes and the leading term 
is $[{\cal M}_{e_{II; 1}, e_{II; 2}, \cdots, e_{II; p}}]_{vir}\cap
c_{p_g}({\cal R}^2\pi_{\ast}{\cal O}_{\cal X})^p$.

\medskip

\noindent Step II: Because both $[{\cal M}_{\underline{C}-\sum e_{II; i}}]_{vir}$
 and $[{\cal M}_{e_{II; 1}, \cdots, e_{II; p}}]_{vir}$ are well defined
 and are independent to $nD$, for the
 whole expression to be $nD$ independent we evaluate the $nD$ independent
 term of 

$$\{c_{a_3-p\cdot q(M)}({\bf U}_{e_{II; 1}, \cdots,
 e_{II; p}}-\oplus_{i\leq p}({\bf R}^0\pi_{\ast}\bigl({
\cal O}_{nD}(nD)\bigr)\oplus {\bf C})
\otimes {\bf H}_{II; i}-{\bf V}'\otimes {\bf H})
\cap [B]_p\}.$$

 Consider the virtual bundle $\omega={\bf U}_{e_{II; 1}, 
e_{II; 2}, \cdots, e_{II; p}}-{\bf V}'\otimes {\bf H}-
\oplus_{i\leq p}({\bf C}\oplus{\bf R}^0\pi_{\ast}\bigl({
\cal O}_{nD}(nD)\bigr))\otimes {\bf H}_{II; i}$. 

The virtual bundle $\omega$ is of virtual rank
$\sum_{i\leq p}\chi({\cal O}_{{\cal X}_b}(e_{II; i}+nD))-
\chi({\cal O}_{\sum {\bf e}_{II; i}|_b}(\underline{C}+nD))-
p(\chi({\cal O}_{{\cal X}_b}(nD))+q)$ (${\cal X}_b$ is the
 fiber of a closed $b\in B$) and by
surface Riemann-Roch formula it is equal to

$$rank(\omega)=-q(M)\cdot p+\sum_{i\leq p}{(e_{II; i})^2}-\underline{C}
\cdot (\sum_{i\leq p}e_{II; i})+\sum_{i<j\leq p}e_{II; i}\cdot e_{II; j}=a_3-p\cdot
 q(M).$$

 The expression $rank(\omega)=a_3-(dim_{\bf C}
\times_B^p {\cal T}_B{\cal X}-dim_{\bf C}B)$ is $nD$ independent!

\begin{defin}\label{defin; noD}
Consider the $nD$ independent term of $c_{rank(\omega)}(\omega)$ and expand it 
as a polynomial of $c_1({\bf H})$, $\sum_{r}\tau_r c_1({\bf H})^r$.
Define the $\tau$ class to be the sum of $\tau_r$, $\tau=\sum_{r\leq rank(\omega)} 
\tau_r$.
\end{defin}

 If we view $c_1({\bf H})$ as a formal variable $z$, then
 $\tau$ can be viewed as $c_{rank(\omega)}(\omega)|_{z=1}^{n=0}$.

\medskip

Recall (e.g. chapter 15, page 281 of [F]) that for a proper morphism $f$ and
 a coherent sheaf ${\cal F}$,  
 $f_{\ast}{\cal F}=\sum_{i\geq 0}(-1)^i{\cal R}^if_{\ast}{\cal F}$.

The following lemma identifies the $nD$ independent term of
$c_{rank(\omega)}(\omega)$.

\begin{lemm}\label{lemm; findout}
The $nD$ independent term of $c_{rank(\omega)}(\omega)$ is equal to
$c_{rank(\omega)}(\oplus_{i\leq p}\pi_{\ast}{\cal O}({\bf e}_{II; i})-
\pi_{\ast}{\cal O}_{\sum_{i\leq p}{\bf e}_{II; i}}(\underline{C})\otimes 
{\bf H})$.
\end{lemm}

\noindent Proof of lemma \ref{lemm; findout}: When $D$ is very ample and
 $n\gg 0$, the Serre vanishing implies 
${\cal R}^0\pi_{\ast}{\cal O}_{\sum {\bf e}_{II; i}}(\underline{C}+nD)
=\pi_{\ast}{\cal O}_{\sum {\bf e}_{II; i}}(\underline{C}+nD)$, etc. This 
enables us to re-express
 $\omega$ as the differences of several direct images. Finally we set $n=0$ in 
 the alternative expression of $\omega$. $\Box$

\bigskip

  Now we may express the $nD$ independent leading term of $g_{a_1, a_2, a_3}$ 
as  

$$[{\cal M}_{\underline{C}-\sum_{i\leq p}e_{II; i}}]_{vir}\times 
\{[{\cal M}_{e_{II; 1}, \cdots, e_{II; p}}]_{vir}
\cap c_{p_g}^p({\cal R}^2\pi_{\ast}{\cal O}_{{\cal X}})\}\cap (\sum_r\tau_r
c_1({\bf H})^r).$$

So the dominating term of the
original localized top Chern class becomes

$$\sum_{r\leq rank(\omega)}\Delta^{!}
\{[{\cal M}_{\underline{C}-\sum_{i\leq p}e_{II; i}}]_{vir}\cap 
c_1({\bf H})^r\times
 [{\cal M}_{e_{II; 1}, \cdots, e_{II; p}}]_{vir}\cap c_{p_g}^p({\bf R}^2\pi_{\ast}
{\cal O}_{\cal X})\cap [B]_p\cap \tau_r\},$$

which is nothing but 
$$\sum_r\{[{\cal M}_{\underline{C}-\sum_{i\leq p}e_{II; i}}]_{vir}\cap 
c_1({\bf H})^r\}\cdot_{id_B}
 \{[{\cal M}_{e_{II; 1}, \cdots, e_{II; p}}]_{vir}\cap c_{p_g}^p
({\bf R}^2\pi_{\ast}
{\cal O}_{\cal X})\cap [B]_p\cap \tau_r\}.$$

 After we cap with $c_1({\bf H})^{dim_{\bf C}B-q+p_g+{\underline{C}^2-
c_1({\bf K}_{{\cal X}/B})\cap \underline{C}\over 2}}$ and 
push forward the top intersection pairing 
to a point $pt$, the top intersection pairing is reduced to
${\cal AFSW}_{{\cal X}\mapsto B}((\times_B^{i\leq p}\pi_i){\ast}\{
[{\cal M}_{e_{II; 1}, \cdots, e_{II; p}}]_{vir}\cap [B]_p\cap \tau\}
\cap c_{p_g}^p({\bf R}^2\pi_{\ast}
{\cal O}_{\cal X}), \underline{C}-\sum_{i\leq p}e_{II; i})$.

\medskip

  By our computation, it clearly obeys guideline 1 on page \pageref{guide}.

\medskip

 To show that it is compatible with the type $I$ theory in [Liu5], [Liu6], we notice that
 $[{\cal M}_{e_{II; 1}, \cdots, e_{II; p}}]_{vir}$ is reduced to
$\cap_{i\leq p}[Y(\Gamma_{e_{k_i}})]=[Y(\Gamma)]$ when the type $II$ classes
 are replaced by $e_{k_i}$, $1\leq i\leq p$. On the other hand,
 we may view $c_1({\bf H})$ as a formal variable $z$, then formally 
$\tau=c_{rank(\omega)}(\omega)|_{z=1}^{n=0}$. When the type $II$ classes
are reduced to the type $I$ classes $e_{k_1}, e_{k_2}, \cdots, e_{k_p}$,
 the argument of theorem 4 of [Liu5] or proposition 18 of [Liu6]
 implies that the effect of 
 $c_{rank(\omega)}(\omega)\cap [B]_p|_{n=0}$ upon the fundamental class
 $[Y(\Gamma)]$ is equal to $c_{rank(\omega)+p\cdot q(M)}({\cal R}^1
\pi_{\ast}{\cal O}_{\sum {\bf e}_{k_i}}(
\underline{C})\otimes {\bf H}-\oplus_i {\cal R}^1\pi_{\ast}
{\cal O}_{{\bf e}_{k_i}}({\bf e}_{k_i}))$ and it is 
equivalent to the top Chern class of an explicit vector bundle
 representative \footnote{By lemma 17 of [Liu6].} of 
$\tau_{\Gamma}\otimes {\bf H}$. So $\tau=
c_{top}(\tau_{\Gamma}\otimes {\bf H})|_{z=1}^{n=0}=c_{total}(\tau_{\Gamma})$.

 The only defect between the degenerated version of type $II$ theory
 and the original type $I$ theory is the expression
$\cap c_{p_g}^p({\cal R}^2\pi_{\ast}{\cal O}_{\cal X})$. This defect roots
 at the discrepancy of their dimension formulae and should be discarded
 when the classes are of type $I$. With this defect removed, the type
$II$ contribution for $\underline{C}-\sum e_{II; i}$ 
is reduced to ${\cal AFSW}_{M_{n+1}\times Y(\Gamma)\mapsto
 M_n\times Y(\Gamma)}(1, C-{\bf M}(E)E-\sum e_{k_i})$ when we take
 $\underline{C}=C-{\bf M}(E)E$.

\medskip

  The object we have identified is independent to $nD$ because of
 the $nD-$independence of $[{\cal M}_{\underline{C}-\sum e_{II; i}}]_{vir}$
 and $[{\cal M}_{e_{II; 1}, \cdots, e_{II; p}}]_{vir}$ and
 the $nD-$independence of $\tau$.

 Thus our construction obeys guidelines 1-3 starting from page \pageref{guide}, 
and the theorem is proved.
 $\Box$

\bigskip

\begin{rem}\label{rem; =0}
When $febd(e_{II; i}, {\cal X}/B)=0$ for all $e_{II; i}$, the
 factor $c_{p_g}^p({\cal R}^2\pi_{\ast}{\cal O}_{\cal X})$ disappears
from the mixed invariant. The mixed invariant identified in the
theorem should be replaced by 

$${\cal AFSW}_{{\cal X}\mapsto B}((\times_B^{i\leq p}
 \pi_i)_{\ast}\bigl([{\cal M}_{e_{II; 1}, \cdots, e_{II; p}}]_{vir}
\cap [B]_p\bigr)\cap \tau, \underline{C}-\sum_{1\leq i\leq p}e_{II; i}).$$
\end{rem}

\begin{rem}\label{rem; tau}
From the above discussion at the end of our main theorem, 
the $\tau$ class defined above is the type $II$ analogue of 
$c_{total}(\tau_{\Gamma})$ defined for type $I$ theory. Its role
 is to balance the rank difference between the family obstruction bundles of
$\underline{C}$ and of $\underline{C}-\sum_{i\leq p}e_{II; i}$. 
The main difference from the type $I$ theory is that for $\tau_{\Gamma}$ we
can represent it as a vector bundle and $c_{top+l}(\tau_{\Gamma})=0$ for all
 $l>0$. But $\omega$ is only a virtual vector bundle of virtual rank
 $-p\cdot q(M)
-\sum_{i\leq p}
\underline{C}\cdot e_{II; i}+\sum_{i\leq p} e_{II; i}^2+\sum_{i<j}e_{II; i}\cdot
 e_{II; j}$. As the type $II$ universal curves ${\bf e}_{II; i}$ may behave
 badly (in the sense described on page \pageref{bizzare}),
 we do not expect to find an explicit bundle representative of $\omega$
generally.
\end{rem}

\subsection{The Remark about the Inductive Scheme of Applying 
Residual Intersection Formula}\label{subsection; inductive}

\bigskip

 At the end of the whole paper, we sketch the extension of the theory to
 an inductive scheme upon a whole hierarchy of collection of type $II$ curves and 
 explain how it fits to our guideline 4 on page \pageref{guide4}.
 As the current discussion is parallel to type $I$ theory in [Liu6], 
we do not intend
to go into the full details here. The reader who wants to get to the
detailed arguments can consult [Liu6]. By combining with the
 main theorem in the current paper, the reader can translate the
 argument to cover the type $II$ case.

\medskip

There are a few reasons why the theory of type $II$ curves has to be
extended beyond a single collection of type $II$ classes.

\medskip

\noindent (1). Exactly parallel to the type $I$ classes, type $II$ curves
can break, or degenerate into a union of irreducible 
curves, while some of them are again
 type $II$ curves. The degenerated configuration will give additional
excess contributions.

\medskip

\noindent (2). The main theorem in the paper 
has been proved under the {\bf special 
assumption} that
 $${\cal M}_{\underline{C}-\sum e_{II; i}}\times_{B'}
{\cal M}_{e_{II; 1}, \cdots, e_{II; p}}\hookrightarrow 
{\cal M}_{\underline{C}}\times_{B'}
{\cal M}_{e_{II; 1}, \cdots, e_{II; p}}$$
 is isomorphic. In section \ref{subsection; stable} this has been interpreted
 equivalently as the injectivity of the bundle map $\pi_{
{\bf P}({\bf V}_{\underline{C}})}^{\ast}{\bf W}_{new}\otimes {\bf H}\mapsto
 \pi_{{\bf P}({\bf V}_{\underline{C}})}^{\ast}
{\bf W}_{\underline{C}}\otimes {\bf H}$ over 
${\cal M}_{\underline{C}-\sum e_{II; i}}\times_{B'}
{\cal M}_{e_{II; 1}, \cdots, e_{II; p}}$.

 In general, the breaking up of the type $II$ classes into
 irreducible components may cause
the inclusion failing to be isomorphic. In other words, a curve in $\underline{C}$
 above the locus of co-existence $\cap_{i\leq p}\pi_i{\cal M}_{e_{II; i}}$
 may not factorize into a curve in $\underline{C}-\sum e_{II; i}$ and
 curves in the sum of $\sum_{i\leq p}e_{II; i}$
 when some curve in $e_{II; i}$ fails to be
 irreducible.

\medskip

 Both theories of type $I$ and type $II$ classes make use of residual
intersection theory of top Chern classes.
 The major difference between type $II$ theory and their type $I$ counter-part
 is that the localized contributions of top Chern classes of type $I$
 classes are identifiable to be mixed family invariants, while
 the localized contribution of top Chern classes of 
type $II$ exceptional classes are typically non-topological. One major
issue we have pointed out 
is that ${\cal M}_{e_{II; 1}, \cdots, e_{II; p}}$ may not always be 
regular and the excess 
invariant contribution within ${\cal M}_{\underline{C}}\times_{B'}
{\cal M}_{e_{II; 1}, \cdots, e_{II; p}}$ may depend on the explicit locus
${\cal M}_{e_{II; 1}, \cdots, e_{II; p}}$ rather than the virtual
 fundamental class $[{\cal M}_{e_{II; 1}, \cdots, e_{II; p}}]_{vir}$.

  Our main theorem on page \pageref{theo; degenerate} demonstrates that despite
 the non-topological nature of the 
localized contribution of top Chern class, 
it can still be expanded algebraically and the dominating term is of the desired
 form: 

 $${\cal AFSW}_{{\cal X}\mapsto B}((\times_B^{i\leq p}\pi_i)_{\ast}\bigl(
[{\cal M}_{e_{II; 1}, \cdots, e_{II; p}}]_{vir}\bigr)\cap 
 c_{p_g}^p({\cal R}^2\pi_{\ast}{\cal O}_{\cal X})\cap [B]_p\cap \tau, 
 \underline{C}-\sum_{i\leq p} e_{II; i}).$$

 This corresponds to the
 invariant contribution proportional
to  $[{\cal M}_{\underline{C}-\sum_{i\leq p}e_{II; i}}]_{vir}\cdot_{id_B}
[{\cal M}_{e_{II; 1}, \cdots, e_{II; p}}]_{vir}$. In this way, the
 non-topological nature of the localized contribution of top Chern classes 
is from the other non-dominating terms from ${\cal M}_{\underline{C}}\times_{B'}
{\cal M}_{e_{II; 1}, \cdots, e_{II; p}}$ away from
 an explicit cycle representative of 
$[{\cal M}_{\underline{C}}]_{vir}\cdot_{id_B}
[{\cal M}_{e_{II; 1}, \cdots, e_{II; p}}]_{vir}$. Unless we impose
additional assumptions, we do not expect
 any vanishing result on these correction terms.

\bigskip

 By re-grouping the correction terms
 with the residual contribution of top Chern class
 $$\int_{\tilde{X}}c_{top}(
\pi_{\tilde{X}}^{\ast}
{\bf W}_{\underline{C}}\otimes {\bf H}\otimes {\cal O}(-{\bf D})\oplus
 {\bf U}_{e_{II; 1}, \cdots, e_{II; p}})\cap [B]_p\cap 
c_1({\bf H})^{dim_{\bf C}B-q+p_g+{\underline{C}^2-c_1({\bf K}_{{\cal X}/B})\cdot
\underline{C}\over 2}},$$

 their total sum is still a topological invariant!

\bigskip

 This interpretation allows us to formulate our scheme inductively.

\medskip

\noindent (i). List all the possible finite collections of type $I$ and 
type $II$ classes satisfying: $\underline{C}\cdot e_{k_i}<0$, 
$e_{k_i}\cdot e_{k_j}\geq 0$, $i\not=j$, 
 $\underline{C}\cdot e_{II; i}<0$, $e_{II; i}\cdot e_{II; j}\geq 0$, $i\not=j$.

These exceptional classes determine exceptional cones in the sense of
 [Liu4].

\medskip

\noindent (ii). Define a partial ordering among all such collections
 based on the inclusions of exceptional cones. 
The partial ordering
  encodes whether one particular collection of exceptional classes
 degenerates into another. In the type $I$ theory, such a partial
ordering has been named $\succ$. Check definition 8 of [Liu6] for details.

\medskip

\noindent (iii). Based on the partial ordering, define a linear ordering among the
 various collections of exceptional classes. This is the analogue of $\models$
 defined in
 [Liu6].

\medskip

\noindent (iv). Blow up ${\bf P}({\bf V}_{\underline{C}})$ 
along the various sub-loci ${\cal M}_{\underline{C}}\times_{B'}
{\cal M}_{e_{k_1}, \cdots, e_{k_p}; e_{II; 1}, \cdots, e_{II; q}}$ by the
reversed linear ordering constructed in (iii). Each time we have to
stabilize the family Kuranishi model of $\underline{C}$ following the
recipe on page \pageref{subsection; stable}.

\medskip

\noindent 
(v). From each localized contribution of top Chern class along

\noindent ${\cal M}_{\underline{C}}\times_{B'}
{\cal M}_{e_{k_1}, \cdots, e_{k_p}; e_{II; 1}, \cdots, e_{II; q}}$, we need to 
 go through
 the computation in theorem \ref{theo; degenerate} and identify its dominating 
term with the appropriated 
mixed family invariant. However, there are two subtle variations here.

\medskip

(v'). The first variation from our main theorem is that under an 
inductive blowing up procedure, the obstruction bundle $\pi_X^{\ast}
{\bf W}_{\underline{C}}\otimes {\bf H}$ has been modified repeatedly.

 As in the type $I$ case, 
we use proposition 9 of [Liu6] to deal with this question.
 The cited proposition demonstrates that the blowing ups performed
 ahead of the given one (under the reversed linear ordering in (iii).) changes the
obstruction bundle in a way which enables us to relate the
 modified bundle with the bundle $
\pi_{{\bf P}({\bf V}_{\underline{C}})}^{\ast}
{\bf W}_{\underline{C}-\sum e_{II; i}}\otimes {\bf H}$, exactly allowing us
 to drop the {\bf special assumption} on page \pageref{special}.

\medskip

A special partial ordering parallel to $\sqsupset$ in definition 15 of [Liu6] 
has to be used to analyze the discrepancy between 
${\cal M}_{\underline{C}-\sum
e_{II; i}}\times_B{\cal M}_{e_{II; 1}, \cdots, e_{II; p}}$ from 
${\cal M}_{\underline{C}}\times_B{\cal M}_{e_{II; 1}, \cdots, e_{II; p}}$.

\medskip

(v''). To avoid over-counting, the inductive blowup process (compare with
 section 5.1. of [Liu6]) incorporates
 the inclusion-exclusion principle (see page 19-20 of [Liu7] for an
 elementary explanation) and we identify the dominating terms of
 the localized contributions to be ``modified'' family invariants.

 For type $I$ classes, their corresponding
 modified family invariants have been defined inductively based on a partial
ordering $\gg$ (definition 11.) in [Liu6]. Please consult definition 13, 14 of 
 the same paper for the definitions of the type $I$ modified family invariants.
For the combinations of type $I$ and
 type $II$ classes, we can generalize $\gg$ and extend the definition of
 modified family invariant accordingly.

 At the end of the inductive procedure, we end up with a modified family
invariant ${\cal AFSW}_{{\cal X}\mapsto B}^{\ast}(1, \underline{C})$.
 It is the subtractions of all the modified family
 invariants attached to the various collections of type $I/II$ exceptional curves 
from the original family invariant 
${\cal AFSW}_{{\cal X}\mapsto B}(1, \underline{C})$.

  This modified family invariant resembles the virtual number count of
 smooth curves in $\underline{C}$ within the given family ${\cal X}\mapsto B$. 

 The above procedure is parallel to the theory of type $I$ exceptional curves.

 If we apply the above scheme to the type $I$ and type $II$ exceptional
 classes of the universal families, we can derive the following result,

\begin{theo}
Given an algebraic surface $M$ and a line bundle $L\mapsto M$. For
 $\delta\leq {c_1^2(L)+c_1({\bf K}_M)\cdot c_1(L)\over 2}+1$, the
 ``virtual number of $\delta$-node nodal curves'' in a generic $\delta$
dimensional linear system of $|L|$ is defined.
\end{theo}

Notice that $\delta!\times$ the above virtual number is interpreted as the
equivalence of smooth curves within the universal family $M_{\delta+1}\mapsto
 M_{\delta}$, represented by some modified family invariant. Unlike the
 universality theorem which works for $5\delta-1$ very ample $L$,
the above result does not guarantee the deformation invariance
 of these virtual numbers of nodal curves.

\bigskip

\subsubsection{The Vanishing Result for $K3$ or $T^4$ and Yau-Zaslow
Formula}\label{subsubsection; 
vanishing}

\medskip

 When we apply the residual intersection theory of type $II$ curves to 
the universal families of $K3$ or $T^4$, we get a surprising vanishing result
 discussed briefly in [Liu7].

Recall that the universality theorem asserts that for a general algebraic
 surface $M$ and a $5\delta-1$ very ample $L$, the $\delta$-node nodal curves
 in a generic $\delta$ dimensional sub-linear system of $|L|$ can 
be expressed in a degree $\delta$ 
universal polynomial of $c_1^2(K_M), c_1(K_M)\cdot c_1(L)$,
$c_1(L)^2$, $c_2(M)$.

  When we consider $M=K3$ (or $T^4$), $K_M$ is trivial and the universal
 polynomial is reduced to a polynomial of $c_1(L)^2$ and $c_2(M)$. On the
 other hand, for rational nodal curves the self-intersection number $c_1(L)^2$
 is constrained by the adjunction formula by $c_1(L)^2=2\delta-2$. 

 So the universal polynomial has been reduced to a degree 
 $\delta$ polynomial of $c_2(M)$ above.

On the other hand, the well known Yau-Zaslow formula [YZ] asserts that when we
 we put all the numbers of $\delta$-node ($\delta\in {\bf N}$) nodal curves 
 into a generating function, it can be identified with the power series 

$$1+\sum_{\delta\in {\bf N}}n_{\delta}q^{\delta}=\{{1\over 
\prod_{i\geq 0} (1-q^i)}\}^{c_2(M)}.$$

The following vanishing result implies that the type $II$ exceptional
 curves within the universal family contribute nothing to the family invariant
 
\noindent ${\cal AFSW}_{M_{n+1}\times \{t_L\}\mapsto M_n\times \{t_L\}}( 1, 
c_1(L)-\sum 2E_i)$.

\begin{theo}\label{theo; zero}
The virtual number of $\delta$-node nodal curves in a linear sub-system of 
$|L|$ 
on an algebraic $K3$ surface is equal to 
${1\over \delta!}{\cal AFSW}_{M_{n+1}\times 
\{t_L\}\mapsto M_n\times \{t_L\}}^{\ast}(1, c_1(L)-2\sum_{i\leq \delta}E_i)$,
 the normalized type $I$ modified algebraic family Seiberg-Witten invariant
 defined in section 5.2., remark 12 of [Liu6].
\end{theo}

\medskip

\noindent Sketch of the 
Proof of theorem \ref{theo; zero}: We apply the above scheme to
both type $I$ and type $II$ exceptional classes on the universal families
 $M_{n+1}\mapsto M_n$. So the $\delta!\times $ the 
virtual number of $\delta-$node nodal curves
 in a linear sub-system of $|L|$ is equal to

${\cal AFSW}_{M_{n+1}\times \{t_L\}\mapsto M_n\times \{t_L\}}( 1, 
c_1(L)-\sum 2E_i)$-correction terms from both type $I$ and type $II$ 
exceptional classes. Each correction term is a modified invariant
 attached to $[{\cal M}_{e_{k_1}, e_{k_2}, \cdots, e_{k_p}; e_{II; 1},
 \cdots, e_{II; p'}}]_{vir}$.

On the other hand, we may distinguish the collections of exceptional classes
 into two subsets. The first subset collects all the collections of type 
$I$ exceptional classes, the second subset collects all the collections
 not entirely of type $I$ exceptional classes. I.e. it contains all the 
collections consisting of either type $II$ exceptional classes or 
consisting of both type $I$ and type $II$ exceptional classes.

 Independent to the details of the $\tau$ class and $Y(\Gamma)$ or 
$[{\cal M}_{e_{II; 1}, \cdots, e_{II; p}}]_{vir}$, by theorem
 \ref{theo; degenerate} the important
characteristic of the mixed invariant involving one or more type $II$ curve
is that there is an additional 
insertion\footnote{Recall that $p_g=1$ for $K3$ surfaces.}
 of $c_1({\cal R}^2\pi_{\ast}{\cal O}_{M_{
\delta+1}})$
 to the family invariant.

 On the other hand, there is the following commutative diagram,

\[
\begin{array}{ccc}
 M_{\delta+1} & \mapsto & M_{\delta}\\
 \Big\downarrow & & \Big\downarrow \\
 M\times M_{\delta} & \mapsto & M_{\delta} 
\end{array}
\]

 and the push-forward of $M_{\delta+1}\mapsto M_{\delta}$ 
factors through $M\times M_{\delta}\mapsto M_{\delta}$.
So ${\cal R}^2\pi_{\ast}{\cal O}_{M_{
\delta+1}}$ can be identified with 
${\cal O}_{M_{\delta}}\otimes H^2(M, {\cal O}_M)$. In particular, this
 implies that $c_1({\cal R}^2\pi_{\ast}{\cal O}_{M_{
\delta+1}})=0$. This implies that all the mixed family invariants involving
 one or more type $II$ classes are identically zero.

 As all the modified invariants are defined inductively by the
 differences of the mixed invariants, all the modified family invariants
 involving type $II$ exceptional curves vanish on the universal families
 of $K3$. Therefore all the correction terms are from the type $I$
exceptional curves. So $\delta!\times$ the virtual number of nodal curves
 collapses to the type $I$ modified family invariant

$${\cal AFSW}_{M_{n+1}\times 
\{t_L\}\mapsto M_n\times \{t_L\}}^{\ast}(1, c_1(L)-2\sum_{i\leq \delta}E_i).$$

The theorem is proved. $\Box$
 
{}

\end{document}